\documentclass[leqno]{amsart}
\usepackage{amsmath,amssymb,mathrsfs,graphicx}
\usepackage{mathptmx,courier}
\usepackage[scaled=0.95]{helvet}

\usepackage{ifpdf}
\ifpdf
\usepackage[pdftex,
            a4paper,
            pdfauthor={Michael Eisermann},
            pdftitle={Arrovian juntas},
            pdfsubject={Classification of arrovian voting system},
            pdfcreator={LaTeX with hyperref},
            pdfproducer={LaTeX with hyperref},
            pdfstartview=FitV,
            bookmarksopen=true,
            bookmarksnumbered,
            colorlinks=true, 
            anchorcolor=black, 
            linkcolor=black, 
            urlcolor=black,
            citecolor=black, 
            pagecolor=black, 
            filecolor=black, 
            menucolor=black, 
            ]{hyperref}
\else
\usepackage[dvips]{hyperref}
\fi

\title{Arrovian juntas}

\author{Michael Eisermann}

\address{Institut Fourier, Universit\'e Grenoble I, France}

\email{Michael.Eisermann@ujf-grenoble.fr}

\urladdr{www-fourier.ujf-grenoble.fr/$^\sim$eiserm}

\date{\today}

\theoremstyle{plain}
  \newtheorem{theorem}{Theorem}
  \newtheorem{lemma}[theorem]{Lemma}
  \newtheorem{proposition}[theorem]{Proposition}
  \newtheorem{corollary}[theorem]{Corollary}
  
\theoremstyle{definition}
  \newtheorem{definition}[theorem]{Definition}
  
  \newtheorem{remark}[theorem]{Remark}
  \newtheorem{example}[theorem]{Example}
  \newtheorem*{notation}{Notation}

\newcommand{\N}{\mathbb{N}}             
\newcommand{\R}{\mathbb{R}}             
\newcommand{\minus}{\smallsetminus}     
\newcommand{\wpref}{\succcurlyeq}       
\newcommand{\spref}{\succ}              
\newcommand{\wferp}{\preccurlyeq}       
\newcommand{\sferp}{\prec}              
\newcommand{\ndiff}{\approx}            
\newcommand{\ncomp}{\perp}              
\newcommand{\PPO}{\mathscr{P}}          
\newcommand{\LPO}{\mathscr{L}}          
\newcommand{\Lex}{\operatorname{Lex}}   
\newcommand{\F}{\mathscr{F}}            
\newcommand{\FF}{\mathbb{F}}            
\renewcommand{\d}{\partial}             
\newcommand{\D}{\mathscr{D}}            

\begin{document} 

\begin{abstract}
  This article explicitly constructs and classifies 
  all arrovian voting systems on three or more alternatives.
  If we demand orderings to be complete, we have, 
  of course, Arrow's classical dictator theorem,
  and a closer look reveals the classification 
  of all such voting systems as dictatorial hierarchies.
  If we leave the traditional realm of complete orderings, the picture changes.
  Here we consider the more general setting where alternatives 
  may be incomparable, that is, we allow orderings that 
  are reflexive and transitive but not necessarily complete.
  Instead of a dictator we exhibit a junta whose internal 
  hierarchy or coalition structure can be surprisingly rich.  
  We give an explicit description of all such voting systems, 
  generalizing and unifying various previous results.
\end{abstract}


\copyrightinfo{2006}{Michael Eisermann}


\maketitle

\vspace{-8mm} \begin{quote} \small \tableofcontents \end{quote} \pagebreak[5]


\section{Introduction and outline of results} \label{sec:Introduction}

\subsection{Motivation and background}

Classifying the objects of an axiomatic theory is a natural 
endeavour whenever it promises to be feasible and meaningful.
In the case of arrovian voting systems, it seems that 
this approach has remained implicit and has not been 
systematically investigated in the published literature.  
This absence is all the more surprising as the arrovian axioms 
were the first to be considered, and characterizations have long 
been accomplished for several other classes of voting rules,
such as simple majority rule \cite{May:1952}
or scoring rules \cite{Smith:1973,Young:1974}.

Ideally, a classification comprises two goals: firstly, 
establish a precise description and compile an exhaustive 
list of all solutions satisfying the requirements, 
and secondly, eliminate possible redundancy by identifying
duplicate descriptions of the same solution.

For complete orderings of at least three alternatives, 
the theorem of K.J.\,Arrow \cite{Arrow:1951} says that every voting system 
satisfying the axioms of unrestricted domain, unanimity, and 
independence of irrelevant alternatives is dictatorial.
As stated, this result does not yet determine the voting system:
if the dictator is indifferent, then all outcomes are still possible.
Nonetheless, the dictator theorem can be used to classify arrovian voting systems: 
by a careful induction argument one can exhibit a dictatorial hierarchy, 
as stated in Theorem \ref{Thm:ArrowRefined}. 

\subsection{From linear to partial orderings}

Arrow's classical result makes essential use of the hypothesis that orderings
be \emph{complete} (also called \emph{linear} or \emph{total}).
While reflexivity and transitivity seem indispensable for social 
orderings, the requirement of completeness is certainly less fundamental.
Driven by Arrow's negative result, it seems worthwhile to drop completeness and 
to consider the more general setting of \emph{partial} orderings 
(also called \emph{quasi-orderings}). 
As we will see, this framework allows for much more flexibility, 
and in particular Arrow's dictator theorem is no longer valid.
It is thus natural to explore the limits 
and to boldly ask: what exactly are the possibilities?

As J.A.\,Weymark \cite{Weymark:1984} pointed out, ``there has been 
surprisingly little work done explicitly on social quasi-orderings''.
He went on to establish an arrovian result by proving the existence 
of an oligarchy, i.e.\ a unique minimal decisive set of voters.
He did not discuss explicit examples, nor did he 
strive for a classification of possible voting systems.

In this article we consider the general setting of partial 
orderings and analyze it as thoroughly as possible.  
We explicitly construct and classify all arrovian voting systems 
on three or more alternatives: we exhibit a junta%
\footnote{ As an aside, in English and many other languages, 
  the word ``junta'' usually has the connotation of military junta, 
  whereas in Spanish the noun ``junta'' can mean a formal assembly 
  in a very general sense.  This spectrum of interpretations corresponds well 
  to the large variety of possible models that we will see later on. }
and precisely describe its internal hierarchy or coalition structure.
To this end we reconsider the notion of (strongly) decisive sets, which 
allows us to classify the simple case (Theorem \ref{Thm:SimpleClassification}).
We then introduce the refined notion of relatively decisive sets,
which allows us to analyze and reconstruct the inner workings 
of every arrovian voting system (Theorem \ref{Thm:GeneralClassification}).

\subsection{Relatively decisive sets} 

We will use fairly standard notation, 
as recalled in \textsection\ref{sec:Notation} below.
Throughout this article, $I$ is the set of individuals (voters), 
and $2^I$ denotes the collection of its subsets.
Up to \textsection\ref{sec:Classification} we will tacitly 
assume $I$ to be finite of size $n$, say $I = \{1,2,\dots,n\}$,
but almost all arguments are valid for infinite sets as well.
We will make this explicit in \textsection\ref{sec:InfiniteSocieties}.

We denote by $A$ the set of alternatives, and we will 
always assume that $A$ contains at least three elements.
The set of all partial preorders on $A$ is denoted by $\PPO_A$.
A voting system is a map $C \colon \PPO_A^I \to \PPO_A$.
We will exclusively be interested in the arrovian case,
where $C$ satisfies the axioms of unanimity and
independence of irrelevant alternatives.

In order to adapt the arrovian approach to partial orderings,
we first revisit the notion of (strongly) decisive subsets of $I$ 
and establish the following result:

\begin{theorem} \label{Thm:SimpleClassification}
  Suppose that the set $A$ contains three or more alternatives and 
  consider a voting system $C \colon \PPO_A^I \to \PPO_A$ that satisfies 
  the axioms of unanimity and independence of irrelevant alternatives.
  Then $C$ also satisfies neutrality and monotonicity.

  Furthermore, there exist unique subsets $K \subseteq J \subseteq I$ such that 
  \begin{enumerate}
  \item
    $a \spref_K b$ implies $a \wpref b$, and $K$ is minimal with respect to this property. \\
    Moreover, each $k \in K$ has veto power:
    if $a \not\wpref_k b$ then $a \not\wpref b$.
  \item
    $a \wpref_J b$ implies $a \wpref b$, and $J$ is minimal with respect to this property. \\
    Moreover, each $j \in J$ has influence:
    if $a \wpref_J b$ and $a \spref_j b$, then $a \spref b$.
  \end{enumerate}

  In other words, $K$ is the minimal decisive set, 
  and $J$ is the minimal strongly decisive set.
  Any pair $(K,J)$ of subsets $\emptyset \ne K \subseteq J \subseteq I$ 
  can be realized in this way, and the trivial case $K = J = \emptyset$ 
  corresponds to the trivial voting system $C(P_1,\dots,P_n) = A \times A$.
  The case $K = J$ corresponds to the voting system
  $C(P_1,\dots,P_n) = \bigcap_{j \in J} P_j$, known as Pareto rule.
\end{theorem}

Knowing the decisive sets does in general 
\emph{not} characterize the voting system,
because the case $K \ne J$ needs further investigation.
In order to obtain complete information, 
we introduce and analyze relatively decisive sets: 
given a subset $N \subseteq I$ we say that $J \subseteq I \minus N$ 
is \emph{decisive relative to $N$} if $a \ndiff_N b$
and $a \spref_J b$ together imply $a \wpref b$,
independently of all other preferences in the profile.
Analogously, $J$ is \emph{strongly decisive relative to $N$} 
if already $a \ndiff_N b$ and $a \wpref_J b$ imply $a \wpref b$.
(For $N = \emptyset$ we recover the absolute notions.)

\subsection{The classification theorem} \label{sub:ClassificationStatement}

We prove that for each $N$ there exists a smallest subset 
$\delta{N} \subseteq I \minus N$ that is decisive relative to $N$.
For ease of notation we also use the notation $\Delta{N} = \delta{N} \sqcup N$, 
which is equivalent to the initial data via $\delta{N} = \Delta{N} \minus N$. 
(We denote by $X \sqcup Y$ the union of two disjoint sets, 
that is, $X \sqcup Y = X \cup Y$ with $X \cap Y = \emptyset$.) 

By definition, $\delta$ and $\Delta$ encode what could be called 
the \emph{coalition structure} underlying the voting system $C$.
This data is then shown to characterize every arrovian voting system,
and we construct several examples in order to illustrate the possibilities.

\begin{theorem} \label{Thm:GeneralClassification}
  Given an arrovian voting system $C \colon \PPO_A^I \to \PPO_A$
  on three or more alternatives, there exists a map $\delta \colon 2^I \to 2^I$
  such that $\delta{N}$ is the minimal decisive set relative to $N$.
  Moreover, $\Delta \colon 2^I \to 2^I$, $\Delta{N} = \delta{N} \sqcup N$ 
  is the unique map enjoying the following properties:
  \begin{enumerate}
  \item
    $N \subseteq \Delta{N}$ and $\Delta{M} \subseteq \Delta{N}$ 
    for all $M \subseteq N \subseteq I$.
    \hfill (Monotonicity)
  \item
    If $\Delta{N} \cap M = \emptyset$ for some $N,M \subseteq I$,
    then $\Delta( N \sqcup M ) = \Delta{N} \sqcup M$.
    \hfill (Minimality)
  \item
    We have $a \wpref b$ if and only if $a \ndiff_N b$ 
    and $a \spref_{\delta{N}} b$ for some $N \subseteq I$.
    \hfill (Decisiveness)
  \end{enumerate}

  Conversely, given an arbitrary map $\Delta \colon 2^I \to 2^I$ 
  satisfying (1), we set $\delta{N} = \Delta{N} \minus N$ and define 
  a map $C_\Delta \colon \PPO_A^I \to \PPO_A$ by the rule (3).
  The resulting voting system $C_\Delta$ satisfies 
  the arrovian axioms of unanimity, neutrality, monotonicity, 
  and independence of irrelevant alternatives.  Condition (2) 
  ensures that $\delta{N}$ is the minimal decisive set relative to $N$.
\end{theorem}

An analogous version holds for infinite societies but its statement 
is necessarily more involved, see Remark \ref{rem:PrincipalVotingSystems}
or more generally Theorem \ref{thm:MeasurableClassification} 
at the end of this article.

The previous theorem establishes a bijective correspondence: 
every arrovian voting system $C \colon \PPO_A^I \to \PPO_A$ is characterized 
by the associated coalition structure $\Delta \colon 2^I \to 2^I$.
Even though conditions (1) and (2) may look technical at first sight, 
they are easy enough to verify for any given map $\Delta$, 
and lend themselves well to the construction of examples: 

\begin{example} \label{exm:NonLexicographic}
  We consider a set $I = \{1,2,3\}$ of three voters and choose 
  $\Delta\emptyset = \{1,2\}$ as the smallest decisive subset.  
  Then we necessarily have $\Delta\{3\} = \{1,2,3\}$, and we also know 
  that $\Delta\{1\}$ and $\Delta\{2\}$ must each contain $\{1,2\}$.
  We choose $\Delta\{1\} = \{1,2,3\}$ and $\Delta\{2\} = \{1,2\}$.  
  This fixes all remaining choices to $\Delta N = \{1,2,3\}$.
  The resulting coalition structure is graphically
  represented in Figure \ref{fig:Coalition}.

  \begin{figure}[hbtp]
    \centering
    \includegraphics{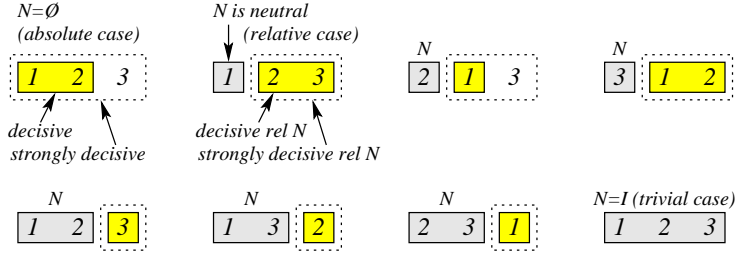}
    \caption{A possible coalition structure for three voters}
    \label{fig:Coalition}
  \end{figure} 

  In words, the associated voting system $C_\Delta$ works as follows.
  Given alternatives $a,b \in A$, voters $1$ and $2$ together 
  can decide either $a \spref b$ or $a \sferp b$ if both agree.
  If $2$ is indifferent, then $1$ alone 
  can decide either $a \spref b$ or $a \sferp b$.
  If $1$ is indifferent, however, then the decision
  is not left to $2$ alone, but to $2$ and $3$:
  they can decide $a \wpref b$ or $a \wferp b$ if both agree.
  In all other cases the conclusion is $a \ncomp b$, 
  that is, $a$ and $b$ are considered incomparable.

  Notice that in this example each voter can influence the outcome:
  we have $a \ndiff b$ if and only if all three voters 
  are indifferent (axiom of strong unanimity).
\end{example}

\begin{remark}
  The apparent complexity is not an artefact of our proof 
  but faithfully reflects the large variety of arrovian voting systems.
  A certain complexity is thus unavoidable: the map $\Delta$ encodes 
  the coalition structure, and, as in game theory or political practice, 
  coalitions are in general a complicated matter.
  Even in small examples such as the previous one, the verbal 
  description can become quite cumbersome, and the formal description 
  using the characteristic map $\Delta$ is usually preferable.
\end{remark}

\subsection{Back to linear orderings} \label{sub:ApplicationToLinearOrderings}

One can weaken the hypotheses by requiring the voting 
system $C$ to be defined only on the set of linear orderings.  
All of our arguments apply with only minor changes,
and we obtain the following result:

\begin{corollary}
  Suppose that $A$ contains three or more alternatives, and 
  let $\LPO_A \subset \PPO_A$ be the set of linear orderings on $A$,
  that is, complete reflexive and transitive relations.
  For every voting system $C \colon \LPO_A^I \to \PPO_A$ satisfying 
  the axioms of unanimity and independence of irrelevant alternatives,
  there exists a unique map $\Delta \colon 2^I \to 2^I$ satisfying 
  the above conditions such that $C = C_\Delta|\LPO_A^I$.
\end{corollary}

The previous corollary can be reformulated more succinctly, without 
explicit reference to the coalition structure $\Delta$, as follows:

\begin{corollary}[Unique arrovian extension] \label{cor:UniqueArrovianExtension}
  Every arrovian voting system $C \colon \LPO_A^I \to \PPO_A$ admits 
  one and only one arrovian extension $\tilde{C} \colon \PPO_A^I \to \PPO_A$ 
  such that $C = \tilde{C}|\LPO_A^I$.
\end{corollary}

This result is rather unexpected.  
Once discovered, it is possible, but not quite obvious, to prove it 
directly, that is, without appealing to the classification. 
It says that we do not gain nor lose any generality 
in choosing the domain $\LPO_A^I$ instead of $\PPO_A^I$.%
\footnote{ We mention in passing that Weymark \cite{Weymark:1984} 
  exclusively considered the case of voting systems $C \colon \LPO_A^I \to \PPO_A$
  where individual preferences are complete while the aggregate preference 
  is allowed to be incomplete, the main advantage being simpler proofs.  
  With hindsight, Corollary \ref{cor:UniqueArrovianExtension} remedies this asymmetry. }

The characterization of arrovian voting systems takes a particularly nice 
and succinct form if we also restrict the range to linear orderings:

\begin{theorem}[Refined version of Arrow's dictator theorem]
 \label{Thm:ArrowRefined}
  Suppose that $A$ contains three or more alternatives.
  Then every map $C \colon \LPO_A^I \to \LPO_A$ that satisfies 
  unanimity and independence of irrelevant alternatives 
  is a lexicographic voting rule according to a unique
  sequence $(k_1,\dots,k_\ell)$ of $\ell \ge 0$ individuals in $I$:
  given any two alternatives, $k_1$ decides on the outcome;
  in case of indifference $k_2$ decides;
  in case of indifference $k_3$ decides, etc.

  Conversely, each sequence $(k_1,\dots,k_\ell)$ defines 
  an arrovian voting system $C \colon \LPO_A^I \to \LPO_A$.
  If we demand that $C$ satisfy the axiom of strong unanimity, 
  then the sequence is necessarily of length $\ell=n$.  There are 
  thus exactly $n!$ voting systems satisfying strong unanimity
  and independence of irrelevant alternatives, 
  each corresponding to a permutation of the set $I$.
\end{theorem}

We obtain this result here as a corollary of the general classification,
but it can also be derived by iterated application of Arrow's 
classical dictator theorem.  Its generalization to infinite societies 
is stated in Theorem \ref{thm:InfiniteDictators}, 
towards the end of this article.

\subsection{How this article is organized}

Section \ref{sec:Notation} recalls the relevant definitions and notation,
in particular concerning orderings (\textsection\ref{sub:Orderings})
and arrovian voting systems (\textsection\ref{sub:ArrovianVotingSystems}).
Section \ref{sec:DecisiveSubsets} reconsiders the notions of decisive
and strongly decisive sets (\textsection\ref{sub:DecisiveSubsets}) and 
establishes their fundamental properties (\textsection\ref{sub:MinimalDecisive}).
This leads to the characterization of juntas without internal structure 
(\textsection\ref{sub:SimpleJuntas}).  In order to treat the general case,
Section \ref{sec:Classification} develops the refined notion of
relatively decisive sets (\textsection\ref{sub:RelativelyDecisive}) 
and discusses lexicographic voting rules as a principal example
(\textsection\ref{sub:Lexicographic}).  These tools allow us to formulate and prove 
the general classification of arrovian juntas (\textsection\ref{sub:Classification})
and its application to the case of linear orderings (\textsection\ref{sub:Applications}).
Section \ref{sec:InfiniteSocieties}, finally, generalizes these results
to infinite societies.


\section{Definitions and notation} \label{sec:Notation}

Our first task is to set the stage so that 
the classification can unfold as smoothly as possible.
While the general axiomatic approach is standard, 
considerable care must be taken to adjust (and motivate) certain details 
of our definitions.  For this reason, and to make this article 
self-contained, we will develop our arguments from scratch.

\subsection{Orderings} \label{sub:Orderings}

In the sequel, $A$ will denote the set of alternatives (decisions, 
proposals, candidates, policies, allocations, issues, etc.)
and preferences will be modelled by orderings on $A$.
More precisely, a \emph{partial preorder} on $A$ is a binary relation 
$P \subseteq A \times A$ that is reflexive and transitive, 
i.e.\ it enjoys the following two properties:
\begin{description}
\item[Reflexivity]
  We have $(a,a) \in P$ for all $a \in A$.
\item[Transitivity]
  Whenever $(a,b) \in P$ and $(b,c) \in P$ then $(a,c) \in P$. 
\end{description}
Throughout this article we will use the term 
\emph{ordering} as synonymous with partial preorder.
We will interpret $(a,b) \in P$ as expressing that 
alternative $a$ is at least as preferable as alternative $b$.
Notice that \emph{indifference} is allowed, that is, $(a,b) \in P$ 
and $(b,a) \in P$ can hold simultaneously for two distinct alternatives $a,b \in A$.  
Moreover, $a$ and $b$ can be \emph{incomparable}, that is, 
neither $(a,b) \in P$ nor $(b,a) \in P$ holds.  

\begin{notation}
  The set of all orderings $P \subseteq A \times A$, 
  i.e.\ partial preorders on $A$, will be denoted by $\PPO_A$.
  We usually write $a \wpref b$ to denote $(a,b) \in P$.
  Indifference ($a \wpref b$ and $b \wpref a$) is denoted $a \ndiff b$.
  Strict preference ($a \wpref b$ but not $b \wpref a$) is denoted 
  $a \spref b$.  We write $a \not\wpref b$ if $(a,b) \notin P$. 
  This includes incomparability (neither $a \wpref b$ nor $b \wpref a$),
  denoted $a \ncomp b$.  
  In summary, we always have $a \ndiff a$, and two distinct alternatives 
  $a \ne b$ are in exactly one of the four relations $a \spref b$ 
  or $b \spref a$ or $a \ndiff b$ or $a \ncomp b$.
\end{notation}

\begin{remark}
  Indifference and incomparability are very different concepts.  
  Notice in particular that indifference is an equivalence relation 
  (i.e.\ reflexive, symmetric, and transitive), whereas incomparability 
  is never reflexive and in general not transitive.

  Incomparability ($a \ncomp b$) is to be interpreted 
  as saying that the comparison between $a$ and $b$ remains 
  undecided on the grounds of the available information.  
  Indifference ($a \ndiff b$), however, expresses the conviction 
  that $a$ and $b$ are equally preferable.  
  Thus, whenever a voting system cannot reach a conclusion 
  between alternatives $a$ and $b$, it seems reasonable 
  to declare them incomparable, rather than equivalent.

  Mathematically speaking, partial preorders are a convenient 
  setting because they allow for sufficient flexibility.
  Admittedly, it is a very general concept and may seem remote from realistic 
  applications.  But then again it nicely models human judgements, 
  where indifference and undecidedness seem to enter quite naturally.  
\end{remark}

\begin{remark}
  Most authors tend to impose stronger conditions 
  by excluding non-comparability, or indifference, or both.  
  In such contexts an ordering $\wpref$ may be required to satisfy:
  \begin{description}
  \item[Completeness]
    Incomparability is excluded: $a \wpref b$ or $b \wpref a$ for all $a,b\in A$.
  \item[Antisymmetry]
    Indifference is excluded: $a \wpref b$ and $b \wpref a$ implies $a = b$. 
  \end{description}
  A complete preorder (as opposed to partial) is also called \emph{linear} 
  or \emph{total}: all alternatives are ranked in a linear fashion, possibly with ties.
  We denote by $\LPO_A \subset \PPO_A$ the set of linear orderings.
  An \emph{order} or \emph{strict ordering} is required to be antisymmetric,
  and a \emph{linear order} is required to be complete \emph{and} antisymmetric.  
  This last notion is the most restrictive one:  it amounts to 
  a strict ranking of all alternatives, without ties, i.e.\ 
  $a_1 \spref a_2 \spref \dots \spref a_m$ in the case of 
  a finite set $A$ of $m$ alternatives.
\end{remark}

\begin{remark}
  Some authors, partly on the basis of empirical psychological evidence, 
  are even willing to sacrifice transitivity, but we will not do so here.
  Transitivity is the fundamental constraint 
  in this whole business and should not be given up.
  All of our arguments will entirely rely on transitivity, and 
  all constructions will carefully preserve this property.
\end{remark}

\subsection{Arrovian voting systems} \label{sub:ArrovianVotingSystems}

Let $\PPO_A$ be the set of all partial preorders on $A$, 
and let $I$ be the set of individuals (voters, agents, committee members).
We will usually assume that $I$ is finite of size $n$.
The cases $n=0$ and $n=1$ are trivial but shall not be excluded.
The subtleties of infinite societies will be discussed
in \textsection\ref{sec:InfiniteSocieties}.

An $n$-tuple $(P_1,\dots,P_n) \in \PPO_A^I$ is called a \emph{preference profile}.  
The \emph{rank aggregation problem} is to single out suitable 
functions $C \colon \PPO_A^I \supseteq \mathscr{D} \to \PPO_A$, 
where the domain $\mathscr{D}$ is some subset of preference profiles,
usually $\mathscr{D} = \PPO_A^I$ or $\mathscr{D} = \LPO_A^I$. 

In voting theory, $C$ is called a \emph{voting system} or \emph{voting rule}
or \emph{social welfare function}: it associates to every profile 
$(P_1,\dots,P_n) \in \mathscr{D}$ of individual preferences $P_1,\dots,P_n$ 
an aggregate preference $P = C(P_1,\dots,P_n)$.
We refer to A.\,Sen \cite{Sen:1986} as a general reference.

\begin{notation}
  It will be convenient to write $a \wpref b$ instead of $(a,b) \in P$,
  and analogously $a \wpref_i b$ as shorthand for $(a,b) \in P_i$.
  The analogous notation $a \spref b$, $a \ndiff b$, $a \ncomp b$, and 
  $a \spref_i b$, $a \ndiff_i b$, $a \ncomp_i b$ with $i \in I$ will also be used.
  Given a subset $J \subseteq I$, we define $a \wpref_J b$ 
  to signify $a \wpref_j b$ for all $j \in J$, and analogously
  $a \spref_J b$ to mean $a \spref_j b$ for all $j \in J$, etc.
\end{notation}

In the sequel we will consider voting systems $C \colon \PPO_A^I \to \PPO_A$, 
which is sometimes called the axiom of \emph{unrestricted domain}, 
because $\mathscr{D} = \PPO_A^I$.  Each such map also satisfies 
what is called the axiom of transitivity, because the range 
is contained in the set of transitive relations.
Among the many desirable properties of such voting systems,
the following axioms have been considered 
as among the most fundamental:

\pagebreak[5]

\begin{description}
\item[Unanimity]
  If all individuals prefer $a$ to $b$
  then so does the aggregate ordering. 
  \begin{description}
  \item[\it Unanimity]
    $a \wpref_I b$ implies $a \wpref b$, i.e.\ 
    if $a \wpref_i b$ for all $i \in I$, then $a \wpref b$.
  \item[\it Strict unanimity]
    $a \wpref_I b$ implies $a \wpref b$, and $a \spref_I b$ implies $a \spref b$.
  \item[\it Strong unanimity]
    $a \wpref_I b$ implies $a \wpref b$.  Moreover, if $a \wpref_I b$ 
    and $a \spref_i b$ for some $i \in I$, then $a \spref b$.
    (Obviously strong unanimity implies strict unanimity.)
\end{description}

\item[Monotonicity]
  If in each individual preference the comparison between $a$ and $b$
  changes in favour of $a$ or remains unchanged, then the same holds 
  true for the aggregate preference.  Formally, consider
  $P = C(P_1,\dots,P_n)$ and $P' = C(P'_1,\dots,P'_n)$:
  if $\{ i \mid a \wpref_i b \} \subseteq \{ i \mid a \wpref'_i b \}$
  and  $\{ i \mid b \wpref_i a \} \supseteq \{ i \mid b \wpref'_i a \}$,
  then $a \wpref b$ implies $a \wpref' b$.

\item[Independence of irrelevant alternatives (IIA)]
  The aggregate ranking between $a$ and $b$ depends 
  only on their individual pairwise rankings.  Formally,
  consider $P = C(P_1,\dots,P_n)$ and $P' = C(P'_1,\dots,P'_n)$:
  if $\{ i \mid a \wpref_i b \} = \{ i \mid a \wpref'_i b \}$
  and  $\{ i \mid b \wpref_i a \} = \{ i \mid b \wpref'_i a \}$,
  then $a \wpref b$ if and only if $a \wpref' b$.
  (This is a consequence of monotonicity.)

\item[Neutrality]
  All alternatives are treated symmetrically, i.e.\ we have
  $C(\rho^* P_1,\dots,\rho^* P_n) = \rho^* C(P_1,\dots,P_n)$
  for all $P_1,\dots,P_n \in \PPO_A$ and 
  every permutation $\rho \colon A \to A$.

\item[Anonymity]
  All individuals are treated symmetrically, i.e.\ we have 
  $C(P_{\sigma 1},\dots,P_{\sigma n}) = C(P_1,\dots,P_n)$
  for all $P_1,\dots,P_n \in \PPO_A$ and 
  every permutation $\sigma \colon I \to I$.
\end{description}

\begin{definition}
  Throughout this article we will suppose that $A$ contains three or 
  more alternatives.  A voting system $C \colon \PPO_A^I \to \PPO_A$
  will be called \emph{arrovian} if it satisfies the axioms 
  of unanimity and independence of irrelevant alternatives.
\end{definition}

\begin{remark}
  As a word of caution, we should point out that we allow 
  the trivial voting system $C_\emptyset \colon \PPO_A^I \to \PPO_A$, 
  $C_\emptyset(P_1,\dots,P_n) = A \times A$, declaring all alternatives 
  equivalent, independently of the preference profile $(P_1,\dots,P_n)$.  
  It satisfies all arrovian axioms except strict unanimity and strong unanimity, 
  but apart from serving as a counter-example it is completely uninteresting.  
  Most authors exclude trivial voting systems by demanding strict unanimity,
  and we will recover this condition as a special case of 
  our classification, see Corollary \ref{cor:StrictAndStrong}.
  The weaker axiom of unanimity is preferable because it behaves well 
  under restriction of the electorate, see Remark \ref{rem:Trivial}.
\end{remark}

\subsection{Immediate consequences} \label{sub:ImmediateConsequences}

We have already noticed that monotonicity implies IIA.
In the presence of three alternatives,
these two axioms become equivalent:

\begin{lemma} \label{lem:NeutralityMonotonicity} 
  If there are at least three alternatives, then 
  unanimity and independence of irrelevant alternatives 
  imply neutrality and monotonicity.
\end{lemma}

\begin{proof}
  Consider a voting system $C \colon \PPO_A^I \to \PPO_A$ 
  satisfying unanimity and independence of irrelevant alternatives. 
  We will first show neutrality, in a stronger form combining both 
  neutrality and independence of irrelevant alternatives:
  \begin{description}
  \item[Strong neutrality]
    Consider $P = C(P_1,\dots,P_n)$ and $P' = C(P'_1,\dots,P'_n)$:
    if $\{ i \mid a \wpref_i b \} = \{ i \mid a' \wpref'_i b' \}$
    and  $\{ i \mid b \wpref_i a \} = \{ i \mid a' \wpref'_i b' \}$,
    then $a \wpref b$ if and only if $a' \wpref' b'$.
  \end{description}
  
  Given two alternatives $a,b \in A$, let $\mathscr{X}(a,b)$ 
  be the set of all pairs $(X,Y)$ such that $a \wpref b$
  with $X = \{ i \mid a \wpref_i b \}$
  and $Y = \{ i \mid a \wferp_i b \}$.
  We claim that $\mathscr{X}(a,b)$ does not depend on 
  the specific alternatives $a$ and $b$, that is, 
  $\mathscr{X}(a,b)$ is the same for all pairs $(a,b)$.

  We will first show that we can replace $b$
  by any other alternative $b' \in A \minus \{a,b\}$,
  according to the following duplicated profile:
  \begin{align*}
    X \cap Y   &: \quad a \ndiff b \ndiff b' \\
    X \minus Y &: \quad a \spref b \ndiff b' \\
    Y \minus X &: \quad a \sferp b \ndiff b' \\
    I \minus ( X \cup Y ) &: \quad a \ncomp b \ndiff b'
  \end{align*}
  Here $b \ndiff b'$ by unanimity,
  hence $a \wpref b \Leftrightarrow a \wpref b'$
  and $a \wferp b \Leftrightarrow a \wferp b'$.
  This proves that $\mathscr{X}(a,b) = \mathscr{X}(a,b')$.
  Analogously we can replace $a$ by any other 
  alternative $a' \in A \minus \{a,b\}$.
  This ultimately leads to $\mathscr{X}(a,b) = \mathscr{X}(a',b')$
  for all pairs $(a,b)$ and $(a',b')$. 

  In order to prove monotonicity, we want to compare two profiles 
  $P = C(P_1,\dots,P_n)$ and $P' = C(P'_1,\dots,P'_n)$ such that 
  \begin{align*}
    X = \{ i \mid a \wpref_i b \} \quad & \subseteq \quad 
    \{ i \mid a \wpref'_i b \} = X' , \\
    Y = \{ i \mid b \wpref_i a \} \quad & \supseteq \quad 
    \{ i \mid b \wpref'_i a \} = Y' .
  \end{align*}
  Assuming $a \wpref b$ we want to show that $a \wpref' b$.
  We choose $b' \in A \minus \{a,b\}$ and construct a third profile 
  $P'' = C(P''_1,\dots,P''_n)$ that duplicates $P$ on the pair $(a,b)$
  and $P'$ on the pair $(a,b')$.  More explicitly, we require 
  $\{ i \mid a \wpref''_i b \} = X$ and $\{ i \mid a \wpref''_i b' \} = X'$,
  as well as $\{ i \mid b \wpref''_i a \} = Y$ and $\{ i \mid b' \wpref''_i a \} = Y'$.
  By our hypothesis $X \subseteq X'$ and $Y \supseteq Y'$ 
  we can also impose $b \wpref''_i b'$ for all $i \in I$.
  Now $a \wpref b$ implies $a \wpref'' b$, by IIA.
  Moreover we have $b \wpref'' b'$ by unanimity, 
  hence $a \wpref'' b'$ by transitivity.
  We conclude that $a \wpref' b$ by IIA and neutrality.
\end{proof}


\section{Decisive subsets} \label{sec:DecisiveSubsets}

\subsection{Definition and first properties} \label{sub:DecisiveSubsets}

Decisive subsets are an important notion for the analysis 
of voting systems.  In the present context two different 
concepts will play a r\^ole:

\begin{definition} \label{def:DecisiveSubset}
  A subset $J \subseteq I$ is called \emph{decisive for the pair $(a,b)$}
  if $a \spref_J b$ implies $a \wpref b$.
  It is called \emph{strongly decisive for the pair $(a,b)$}
  if already $a \wpref_J b$ implies $a \wpref b$.
\end{definition}

Notice that for $J$ to be decisive we only demand that $a \spref_J b$ 
imply $a \wpref b$.  The expected stronger conclusion $a \spref b$
will be a consequence, shown in Proposition \ref{prop:VetoPower},
provided that the voting system is non-trivial.%
\footnote{ For the trivial voting system, the empty set 
  $J = \emptyset$ is decisive: $a \spref_J b$ is always true 
  and implies $a \wpref b$. (But not $a \spref b$!)  
  Conversely, if $J = \emptyset$ is decisive, 
  then the voting system is necessarily trivial.  
  The trivial case is thus completely harmless, 
  and there is no need to disallow it.  
  See Remark \ref{rem:Trivial} for an argument 
  in favour of this choice. }

\begin{remark} \label{rem:WorstCase}
  The axiom of unanimity states that $I$ is 
  strongly decisive for every pair $(a,b)$.
  Obviously strongly decisive implies decisive.
  Given neutrality, if a subset is (strongly) decisive for \emph{one} 
  pair $(a,b)$, then it is (strongly) decisive for all pairs.
  Given monotonicity, it suffices to consider the worst case:
  \begin{enumerate}
  \item \label{item:WorstCaseDecisive}
    A subset $J \subseteq I$ is decisive for $(a,b)$ 
    if and only if $a \spref_J b$ on $J$ and $a \sferp_{\bar{J}} b$ 
    on the complement $\bar{J} = I \minus J$ imply $a \wpref b$.
  \item \label{item:WorstCaseStronglyDecisive}
    A subset $J \subseteq I$ is strongly decisive for $(a,b)$ 
    if and only if $a \ndiff_J b$ on $J$ and $a \sferp_{\bar{J}} b$ 
    on the complement $\bar{J} = I \minus J$ imply $a \wpref b$.
  \end{enumerate}
\end{remark}

\begin{proposition}[Exclusion property] \label{prop:ExclusionProperty}
  If $J \subseteq I$ is strongly decisive,
  then the members of the complement $I \minus J$
  have no influence whatsoever on the outcome.
\end{proposition}

\begin{proof}
  Given two alternatives $a,b \in A$, we will show 
  that the outcome for $a$ and $b$ depends only on 
  the sets $X = \{j \in J \mid a \wpref_j b \}$
  and  $Y = \{j \in J \mid a \wferp_j b \}$.
  Consider a third alternative $a' \in A \minus \{a,b\}$ and assume 
  $a' \wpref_X b$ and $a' \wferp_Y b$, as well as $a \ndiff_J a'$.  
  Then $a \ndiff a'$ because $J$ is strongly decisive.  
  The preferences on $\{a',b\}$ determine the outcome on $\{a',b\}$,
  and thus the outcome on $\{a,b\}$ because $a \ndiff a'$.
  We can thus arbitrarily modify the preferences of $I \minus J$ 
  on $\{a,b\}$ without changing the outcome.
\end{proof}

\begin{proposition}[Intersection property] \label{prop:IntersectionProperty}
  Suppose that the set $A$ contains at least three alternatives and 
  that the voting system $C \colon \PPO_A^I \to \PPO_A$ satisfies 
  unanimity and independence of irrelevant alternatives.
  If $J_1$ and $J_2$ are two (strongly) decisive subsets of $I$,
  then so is their intersection $J = J_1 \cap J_2$.
  If the set $I$ of voters is finite, then the intersection of all 
  (strongly) decisive subsets is the unique minimal (strongly) decisive subset.
\end{proposition}

\begin{proof}
  Assume that $J_1$ and $J_2$ are decisive.
  Given $a \spref_J b$ we want to show that $a \wpref b$.
  We can choose a third alternative $x \in A \minus \{a,b\}$
  and arrange $a \spref_{J_1} x$ and $x \spref_{J_2} b$ as follows:
  \begin{align*}
    J            &: \quad a \spref x \spref b \\
    J_1 \minus J &: \quad a,b \spref x \\
    I \minus J_1 &: \quad x \spref a,b
  \end{align*}
  We obtain $a \wpref x$ because $J_1$ is decisive,
  and $x \wpref b$ because $J_2$ is decisive,
  hence $a \wpref b$ by transitivity.  
  This proves that $J = J_1 \cap J_2$ is 
  decisive for $(a,b)$, hence decisive by neutrality.
  The argument for strongly decisive sets is analogous,
  by replacing ``$\spref$'' with ``$\wpref$''.
\end{proof}

\subsection{Characterizing minimality} \label{sub:MinimalDecisive}

We begin with a sufficient criterion for decisiveness:

\begin{lemma} \label{lem:DecisiveSupport}
  If $a \wpref b$ then the supporting set
  $K = \{ i \in I \mid a \wpref_i b \}$ is decisive.
\end{lemma}

\begin{proof}
  Consider $J = \{ i \in I \mid a \sferp_i b \}$
  and $L = \{ i \in I \mid a \ncomp_i b \}$.
  We can choose a third alternative $c \in A \minus \{a,b\}$
  and assume $b \spref_K c$ and $b \sferp_J c$ and $b \ncomp_L c$.
  Transitivity entails $a \spref_K c$ and $a \sferp_J c$.
  On $L$, however, there is no restriction on the relation 
  between $a$ and $c$, and so we can choose $a \sferp_L c$.
  \begin{align*}
    K &: \quad a \wpref b \spref c \\
    J &: \quad a \sferp b \sferp c \\
    L &: \quad a \ncomp b \ncomp c, \quad a \sferp c
  \end{align*}
  By neutrality and monotonicity, $a \wpref b$ implies 
  $b \wpref c$, hence $a \wpref c$ by transitivity.  
  Now Remark \ref{rem:WorstCase}\,\eqref{item:WorstCaseDecisive}
  applies and we conclude that $K$ is decisive for $(a,c)$, hence decisive.
\end{proof}


\begin{proposition} \label{prop:VetoPower}
  For a decisive set $K \subseteq I$ the following conditions are equivalent:
  \begin{enumerate}
  \item $K$ is minimal, that is, 
    if $K' \subseteq K$ is decisive then $K' = K$.
  \item $K$ is the smallest decisive set, 
    that is, if $K'$ is decisive then $K \subseteq K'$.
  \item Each member $k \in K$ has veto power
    in the sense that $a \not\wpref_k b$ implies $a \not\wpref b$.
  \end{enumerate}
  As a consequence, $a \spref_K b$ implies $a \spref b$,
  provided that the voting system is non-trivial.
\end{proposition}

\begin{proof}
  The implication $(1) \Rightarrow (2)$ follows
  from Proposition \ref{prop:IntersectionProperty}:
  if $K'$ is decisive, then $K \cap K'$ is decisive,
  hence $K = K \cap K'$ by minimality of $K$.
  This means that $K \subseteq K'$.

  The implication $(2) \Rightarrow (3)$ follows
  from Lemma \ref{lem:DecisiveSupport}:  if $a \wpref b$ 
  then $\{ i \in I \mid a \wpref_i b \}$ is decisive.
  It necessarily contains the smallest decisive subset $K$, 
  thus $a \wpref_K b$.  By contraposition, 
  if $a \not\wpref_k b$ for some $k \in K$ then $a \not\wpref b$.

  The implication $(3) \Rightarrow (1)$ is clear:
  if $K' \subseteq K$ is decisive, then $k \in K \minus K'$ has no veto power.
  This being excluded, we necessarily have $K' = K$, i.e.\ $K$ is minimal.
\end{proof}

\begin{remark}[Oligarchies]
  In the terminology of Weymark \cite{Weymark:1984}, following Gibbard \cite{Gibbard:1969},
  the minimal decisive set $K \subseteq I$ is an \emph{oligarchy}, in the sense that 
  its members can impose $a \spref b$ in the case of unanimity $a \spref_K b$, while each individual 
  member $k \in K$ has veto power, i.e.\ $a \sferp_k b$ implies $a \not\wpref b$.
  Weymark studied only voting systems $C \colon \LPO_A^I \to \PPO_A$
  and thus excluded incomparabilities in the individual preferences.
  The preceding proposition proves the stronger version
  that $a \not\wpref_k b$ implies $a \not\wpref b$, which 
  covers both $a \sferp_k b$ and $a \ncomp_k b$.
\end{remark}

\begin{remark}[Arrow's theorem]
  Our arguments also hold for voting systems $C \colon \LPO_A^I \to \PPO_A$
  defined on the subset $\LPO_A \subset \PPO_A$ of linear orderings.
  Restricting to this domain simplifies the proofs, but does not 
  change the results in an essential way.  

  Restricting the range, however, is more severe.
  If we consider $C \colon \LPO_A^I \to \LPO_A$, then the minimal 
  decisive subset $K$ cannot contain more than one individual: 
  if two individuals $i,j \in K$ have veto power, then 
  $a \spref_i b$ and $a \sferp_j b$ would imply $a \ncomp b$.
  We thus obtain, as a by-product, Arrow's classical dictator theorem.
  (We will come back to this in \textsection\ref{sub:Applications}.)
\end{remark}

Analogously, we have a sufficient criterion for strong decisiveness:

\begin{lemma} \label{lem:StronglyDecisiveSupport}
  If $a \ndiff b$ then the supporting set
  $N = \{ i \in I \mid a \ndiff_i b \}$ is strongly decisive.
\end{lemma}

\begin{proof}
  Consider $J = \{ i \in I \mid a \spref_i b \}$
  and $K = \{ i \in I \mid a \sferp_i b \}$
  and $L = \{ i \in I \mid a \ncomp_i b \}$.
  We choose a third alternative $a' \in A \minus \{a,b\}$ 
  and consider the following profile:
  \begin{align*}
    N &: \quad a' \ndiff a \ndiff b \\
    J &: \quad a' \spref a \spref b \\
    K &: \quad a \sferp a' \sferp b \\
    L &: \quad a \sferp a', \quad a \ncomp b, \quad a' \ncomp b
  \end{align*}
  By neutrality, $a \ndiff b$ implies $a' \ndiff b$, whence $a \ndiff a'$.
  Now Remark \ref{rem:WorstCase}\,\eqref{item:WorstCaseStronglyDecisive}
  applies and we conclude that $N$ is strongly decisive for $(a,a')$, 
  hence strongly decisive.
\end{proof}

\begin{proposition} \label{prop:MinimalStronglyDecisive}
  For a strongly decisive set $J \subseteq I$ the following conditions are equivalent:
  \begin{enumerate}
  \item $J$ is minimal, that is, 
    if $J' \subseteq J$ is strongly decisive then $J' = J$.
  \item $J$ is the smallest strongly decisive set, 
    that is, if $J'$ is strongly decisive then $J \subseteq J'$.
  \item Each member $j \in J$ has influence in the sense that 
    $a \wpref_J b$ and $a \spref_j b$ imply $a \spref b$.
  \end{enumerate}
  In particular, the voting system satisfies strong unanimity if and only if 
  the set $I$ of all voters is the minimal strongly decisive subset.
\end{proposition}

\begin{proof}
  The implication $(1) \Rightarrow (2)$ follows
  from Proposition \ref{prop:IntersectionProperty}:
  if $J'$ is strongly decisive, then so is $J \cap J'$,
  hence $J = J \cap J'$ by minimality of $J$.
  This means that $J \subseteq J'$.

  The implication $(2) \Rightarrow (3)$ follows
  from Lemma \ref{lem:StronglyDecisiveSupport}: if $a \ndiff b$ 
  then $\{ i \in I \mid a \ndiff_i b \}$ is strongly decisive.
  It necessarily contains the smallest strongly decisive subset $J$, 
  thus $a \ndiff_J b$.  By contraposition, if $a \not\ndiff_j b$ 
  for some $j \in J$ then $a \not\ndiff b$.  We conclude that
  $a \wpref_J b$ and $a \spref_j b$ for some $j \in J$ imply $a \spref b$.

  The implication $(3) \Rightarrow (1)$ is clear:
  if $J' \subseteq J$ is strongly decisive, 
  then $a \ndiff_{J'} b$ implies $a \ndiff b$,
  and $j \in J \minus J'$ has no influence on the outcome.
\end{proof}

\begin{corollary} \label{cor:StrictAndStrong}
  For an arrovian voting system $C \colon \PPO_A^I \to \PPO_A$ 
  on three or more alternatives, let $K \subseteq J \subseteq I$ be 
  the minimal decisive and strongly decisive subsets, respectively.  
  \begin{enumerate}
  \item
    $C$ satisfies the axiom of strict unanimity if and only if 
    $K \ne \emptyset$, i.e.\ $C$ is non-trivial.
  \item
    $C$ satisfies the axiom of strong unanimity if and only if $J = I$.
    \qed
  \end{enumerate}
\end{corollary}

\begin{remark}[Revelation principle]
  From the very beginning of our discussion we have supposed that 
  the individual preferences $P_1,\dots,P_n \in \PPO_A$ are known, 
  tacitly assuming that all individuals truthfully reveal their preferences.  
  From a game-theoretic point of view, this assumption is consistent 
  if revealing the true preferences is a Nash equilibrium: 
  no individual is better off by declaring another preference,  
  or in other words, no individual has an incentive to lie.
  Moreover, individuals have an incentive to truthfully reveal their 
  preferences if each individual preference can potentially decide the outcome.
  For this strong form of Nash equilibrium it is necessary 
  that $C$ be strongly unanimous: if $J \subsetneq I$ is strongly decisive, 
  then some individuals are systematically excluded from the decision 
  process and have no incentive whatsoever to reveal their preferences.
\end{remark}

\subsection{Juntas without internal structure} \label{sub:SimpleJuntas}

The preceding notions being in place, we can now begin
to classify the precise structure of arrovian voting systems.
We will first treat the simple case where decisive 
and strongly decisive subsets co\"incide.
The general case is more complex and 
will be treated in the next section.

\begin{example}
  Choose a subset $J \subseteq I$ and define 
  the voting system $C_J \colon \PPO_A^I \to \PPO_A$ 
  by $C_J(P_1,\dots,P_n) := \bigcap_{j \in J} P_j$.
  This implements the \emph{Pareto rule} based on $J$, which means that 
  $a \wpref b$ if and only if $a \wpref_j b$ for all $j \in J$.
  This can be seen as the ``greatest common ordering''
  unanimously agreed upon by all members of $J$.  
  In the trivial case $J = \emptyset$ we recover the trivial 
  voting system $C_\emptyset(P_1,\dots,P_n) = A \times A$.
\end{example} 

It is easily verified that $C_J$ satisfies the axioms of unanimity, 
neutrality, monotonicity, and independence of irrelevant alternatives.
Notice also that $J$ is the smallest decisive subset, and 
at the same time the smallest strongly decisive subset.
This property is characteristic in the following sense:

\begin{proposition}[Classification of juntas without internal structure] \label{prop:SimpleJuntas}
  Suppose that the set $A$ contains at least three alternatives and 
  consider a voting system $C \colon \PPO_A^I \to \PPO_A$ that satisfies 
  the axioms of unanimity and independence of irrelevant alternatives.
  Let $K \subseteq I$ be the smallest decisive subset, and 
  let $J \subseteq I$ be the smallest strongly decisive subset.
  Then we have the double inclusion $C_J \subseteq C \subseteq C_K$, that is,
  \[
  \textstyle
  \bigcap_{j \in J} P_j 
  \quad\subseteq\quad
  C(P_1,\dots,P_n) 
  \quad\subseteq\quad 
  \bigcap_{k \in K} P_k
  \]
  for all preference profiles $(P_1,\dots,P_n) \in \PPO_A^I$.
  Equality $C = C_J$ or $C = C_K$ holds if and only if $K = J$, 
  that is, decisive and strongly decisive subsets co\"incide.
\end{proposition}

\begin{proof}
  On the one hand, the fact that $J$ is strongly decisive is equivalent 
  to the inclusion $\bigcap_{j \in J} P_j \subseteq C(P_1,\dots,P_n)$
  for all profiles $(P_1,\dots,P_n) \in \PPO_A^I$.
  On the other hand, according to Proposition \ref{prop:VetoPower},
  each member $k \in K$ has veto power in the sense that 
  $a \not\wpref_k b$ implies $a \not\wpref b$.  This is equivalent 
  to $C(P_1,\dots,P_n) \subseteq \bigcap_{k \in K} P_k$.
  Obviously, if $C = C_J$ or $C = C_K$ then $K = J$.  Conversely, 
  if $K = J$, then the above inclusions show that $C = C_J = C_K$.
\end{proof}

\begin{proposition} \label{prop:RestrictionToLinearOrderings}
  The statement of the preceding proposition holds verbatim 
  for voting system $C \colon \LPO_A^I \to \PPO_A$ defined 
  only on the subset $\LPO_A \subset \PPO_A$ of linear orderings.
\end{proposition}

\begin{proof}
  We claim that the arguments developed for partial orderings 
  still apply to linear orderings, i.e.\ all the necessary constructions 
  can be carried out within $\LPO_A^I \subset \PPO_A^I$.
  
  First of all, the definition of (strongly) decisive sets applies 
  to linear orderings (Definition \ref{def:DecisiveSubset}).  
  The intersection property remains valid and ensures the existence 
  of a unique minimal decisive set (Proposition \ref{prop:IntersectionProperty}).
  We still have the characterization of the minimal decisive set
  in terms of veto power (Proposition \ref{prop:VetoPower});
  the proof is somewhat simplified by the stronger hypothesis that 
  individual preferences no longer present incomparabilities.
\end{proof}

\begin{corollary}
  The only arrovian voting systems $\PPO_A^I \to \PPO_A$
  (or $\LPO_A^I \to \PPO_A$, respectively) on three or more 
  alternatives that satisfy anonymity are the trivial voting 
  system $C_\emptyset$ and the strong Pareto rule $C_I$.
  Only $C_I$ satisfies anonymity and strong unanimity.
\end{corollary}

\begin{proof}
  Let $C \colon \PPO_A^I \to \PPO_A$ or $C \colon \LPO_A^I \to \PPO_A$
  be an arrovian voting system and let $K$ be its minimal decisive set.
  Anonymity implies $K = \emptyset$ or $K = I$, thus $C = C_\emptyset$ or $C = C_I$.
\end{proof}

\begin{remark}
  Although $C_J$ fulfils the arrovian requirements, 
  it is not at all satisfactory:
  \begin{itemize}
  \item
    In the extreme case $J = \{j\}$ the voting system $C_J$
    is dictatorial, with the individual $j$ as dictator.
    In order to integrate more individuals in the decision process,
    we have to enlarge the set $J$, but there is a price to pay:
    we can no longer ensure a complete ordering of alternatives.
    Even if we start out with complete orderings $P_j$, their 
    intersection $P = \bigcap_{j \in J} P_j $ may not be complete.  
    In the worst case, we obtain the trivial ordering 
    $P = A \times A$. 
  \item
    In the other extreme, $J=I$, the map $C_I$ is the \emph{strong Pareto rule}.
    It is the consensus in the strongest possible sense.
    Since each individual $i \in I$ has veto power, a single 
    deviant voter can annihilate this over-ambitious consensus: 
    by insisting on $a \ncomp_i b$ each one can enforce the outcome $a \ncomp b$.
    Even if we exclude incomparability and consider $C_I \colon \LPO_A^I \to
    \PPO_A$, the same happens for two opposing voters $i,j$ 
    with preferences $a \spref_i b$ and $a \sferp_j b$.
  \end{itemize}
  In particular we see that the voting system $C_J$ wastes 
  a lot of information: if $J$ is small, many individuals are excluded;
  if $J$ is large, then the conclusions are diluted, to the extent 
  that quite often no conclusion whatsoever can be drawn.
  (Citizens of the European Union can testify recent examples.)
  This is the classical dilemma between equity and efficiency,
  which is the core problem in voting theory. 
\end{remark}

\begin{remark}
  The general setting of partial orderings highlights the r\^ole of
  completeness in Arrow's theorem: while in general any subset 
  $K \subseteq I$ can be an oligarchy, completeness forces oligarchies 
  to shrink to a single individual, whence the dictator theorem. 
\end{remark}


\section{Classification of arrovian juntas} \label{sec:Classification}

\subsection{Relatively decisive subsets} \label{sub:RelativelyDecisive}

The preceding Proposition \ref{prop:SimpleJuntas} characterizes arrovian 
voting systems in which decisive and strongly decisive subsets co\"incide.
In general they differ, as shown by Example \ref{exm:Lexicographic} below.
In order to classify all possibilities, a more detailed analysis is thus required.
The appropriate refinement is that of relatively decisive subsets:

\begin{definition} \label{def:RelativelyDecisiveSubset}
  Let $N \subseteq I$ be a subset and consider two alternatives $a,b \in A$.  
  A subset $J \subseteq I \minus N$ is \emph{decisive relative to $N$} 
  for the pair $(a,b)$ if $a \ndiff_N b$ and $a \spref_J b$ together 
  imply $a \wpref b$, independently of all other preferences.
  It is \emph{strongly decisive relative to $N$} for the pair $(a,b)$ 
  if $a \ndiff_N b$ and $a \wpref_J b$ together imply $a \wpref b$.
\end{definition}

Given neutrality, if a subset $J$ is (strongly) decisive relative to $N$ 
for \emph{one} pair $(a,b)$, then it is (strongly) decisive 
relative to $N$ for all pairs.  In the latter case $J$ 
is simply called \emph{(strongly) decisive relative to $N$},
without reference to any pair $(a,b)$.
Notice that in the case $N = \emptyset$ we recover 
the absolute version of (strong) decisiveness
as in Definition \ref{def:DecisiveSubset}.

The interpretation of relative decisiveness is as follows:
if the members of $N$ declare themselves neutral in the sense 
that they regard $a$ and $b$ as being equivalent, then
the decision is left to the members of the complement $I \minus N$.
This can be formalized as follows: the inclusion $I \minus N \hookrightarrow I$ 
induces an inclusion $\phi^N \colon \PPO_A^{I \minus N} \hookrightarrow \PPO_A^I$, 
$(P_i)_{i \in I \minus N} \mapsto (P_i)_{i \in I}$ 
by extending with $P_i = A \times A$ for all $i \in N$.  
We can thus define the restricted voting system 
$C^N := C \circ \phi^N \colon \PPO_A^{I \minus N} \hookrightarrow \PPO_A^I \to \PPO_A$.

\begin{proposition}[Restriction of the electorate] \label{prop:RestrictionOfElectorate}
  All arrovian axioms (except strict unanimity) are hereditary in the sense 
  that they remain valid when passing from $C$ to the restriction $C^N$.
  A subset $J \subseteq I \minus N$ is (strongly) decisive relative to $N$
  if and only if $J$ is (strongly) decisive for the voting system $C^N$.
  All properties established for (strongly) decisive sets thus 
  carry over to (strongly) decisive sets relative to $N$.
  In particular, if $I$ is finite, then for each $N \subseteq I$ 
  there exists a unique minimal decisive subset relative to $N$, 
  denoted $\delta{N} \subseteq I \minus N$.
  \qed
\end{proposition}

For ease of notation it is sometimes more convenient to work with
$\Delta{N} = \delta{N} \sqcup N$.  Since $N \cap \delta{N} = \emptyset$, 
we can recover the initial data via $\delta{N} = \Delta{N} \minus N$.

\begin{example}[Pareto rule]
  Consider $C_J \colon \PPO_A^I \to \PPO_A$,
  $C_J(P_1,\dots,P_n) := \bigcap_{j \in J} P_j$
  for some fixed subset $J \subseteq I$.  
  Here we find $\delta{N} = J \minus N$ and 
  $\Delta{N} = J \cup N$ for all $N \subseteq I$.
\end{example}

\begin{remark} \label{rem:Trivial}
  In order to restrict the electorate as explained in 
  Proposition \ref{prop:RestrictionOfElectorate}, 
  we have to allow trivial voting systems:
  even if $C$ is non-trivial, it may well be that $C^N$ is trivial.
  This is the case if and only if $N$ is strongly decisive for $C$,
  see Proposition \ref{prop:ExclusionProperty}.  In particular, the axiom 
  of strict unanimity is not hereditary and thus technically quite cumbersome.

  Trivial voting systems could be avoided by demanding strong unanimity.  
  This axiom is hereditary and our approach could thus be based on this stronger condition.  
  For a finite society $I$ both choices are essentially equivalent:
  each voting system $\PPO_A^I \to \PPO_A$ satisfying unanimity 
  can be restricted to $\PPO_A^J \to \PPO_A$ satisfying strong unanimity, 
  where $J \subseteq I$ is the minimal strongly decisive subset.
  (See Propositions \ref{prop:ExclusionProperty} and \ref {prop:MinimalStronglyDecisive}.)
\end{remark}

\subsection{Lexicographic voting rules} \label{sub:Lexicographic}

The voting rule $C_J$ can be refined
by introducing extra structure.  The idea is,
as could be suspected, that all junta members 
are equal, but some are more equal than others:

\begin{example}[Lexicographic voting rule] \label{exm:Lexicographic}
  Let $\Omega = \{ J_1 \subset \dots \subset J_\ell \}$
  be an ascending chain of non-empty subsets $J_\lambda \subseteq I$.
  We can then define $\Lex_\Omega \colon \PPO_A^I \to \PPO_A$ as follows: 
  Given two alternatives $a,b \in A$, we set $a \ndiff b$ 
  if and only if $a \ndiff_{J} b$ for all $J \in \Omega$.
  Otherwise, let $J \in \Omega$ be the smallest set
  for which we do not have $a \ndiff_J b$; we then set 
  $a \wpref b$ if and only if $a \wpref_J b$, and symmetrically 
  $a \wferp b$ if and only if $a \wferp_J b$.
  This is the \emph{lexicographic} voting rule associated to the chain $\Omega$:
  one can interpret $\Omega$ as a hierarchically ordered junta, 
  with $J_1$ being the inner circle, $J_2$ being the enlarged
  inner circle, etc.  Only in the case of indifference 
  is the decision passed down in the hierarchy.
\end{example}

\begin{example}[Strong lexicographic voting rule] \label{exm:StrongLexicographic}
  Given a chain $\Omega = \{ J_1 \subset \dots \subset J_\ell \}$
  we set $J_0 = \emptyset$ and $J_\lambda = J_\ell$ for all $\lambda \ge \ell$.
  The lexicographic voting rule $\Lex_\Omega \colon \PPO_A^I \to \PPO_A$ 
  can then be rephrased as follows: we have $a \wpref b$ if and only if 
  there exists $J = J_\lambda$ and $K = J_{\lambda+1} \minus J_\lambda$ 
  such that $a \wpref_J b$ and $a \wpref_K b$ but not $a \ndiff_K b$.  
  By a slight modification, we define another voting rule
  $\Lex'_\Omega \colon \PPO_A^I \to \PPO_A$ as follows: 
  we set $a \wpref b$ if and only if there exists $J = J_\lambda$
  and $K = J_{\lambda+1} \minus J_\lambda$ such that $a \wpref_J b$ 
  and $a \spref_K b$.  This is called the \emph{strong} 
  lexicographic voting rule associated to the chain $\Omega$.
  As long as $a \wpref_J b$, the decision is passed down in the hierarchy
  until the relatively decisive set $K$ agrees on a strict preference $a \spref_K b$.
\end{example}

\begin{proposition} \label{prop:LexicographicRule}
  For every chain $\Omega$ 
  the voting rules $\Lex_\Omega$ and $\Lex'_\Omega$ satisfy the axioms of 
  unanimity, neutrality, monotonicity, and independence of irrelevant alternatives.

  For $\Lex_\Omega$ we find $\delta_\Omega{N} = \emptyset$ if $J_\ell \subseteq N$,
  and otherwise $\delta_\Omega{N} = J_\lambda \minus N$ where $J_\lambda \in \Omega$ 
  is the smallest subset such that $J_\lambda \not\subseteq N$.
  For $\Lex'_\Omega$ we find $\delta'_\Omega{N} = \emptyset$ if $J_\ell \subseteq N$, 
  and otherwise $\delta'_\Omega{N} = J_{\lambda+1} \minus N$ where $J_\lambda \in \Omega$ 
  is the smallest subset such that $J_{\lambda+1} \cap N \subseteq J_\lambda$.

  In both cases we obtain $\Delta\emptyset = J_1$, $\Delta{J_1} = J_2$, 
  \dots, $\Delta{J_{\ell-1}} = J_\ell$, $\Delta{J_\ell} = J_\ell$,
  in particular the voting rule $\Lex_\Omega$ or $\Lex'_\Omega$ 
  uniquely determines the chain $\Omega$ used in its construction.
  \qed
\end{proposition}

More generally, we can associate to each arrovian voting system 
$C \colon \PPO_A^I \to \PPO_A$ a chain $\Omega$ of subsets of $I$
by setting $J_0 := \emptyset$ and inductively $J_{\lambda+1} := \Delta{J_\lambda}$. 
This data partitions arrovian voting systems into certain families,
each associated to some chain $\Omega$.  In each such family, 
$\Lex'_\Omega$ and $\Lex_\Omega$ are the least and 
the greatest element, respectively:

\begin{proposition} \label{prop:LexicographicInclusion}
  For every arrovian voting system $C \colon \PPO_A^I \to \PPO_A$ 
  there exists a unique chain $\Omega$ of subsets of $I$
  such that $\Lex'_\Omega \subseteq C \subseteq \Lex_\Omega$,
  namely the chain inductively defined by $J_0 := \emptyset$ 
  and $J_{\lambda+1} := \Delta{J_\lambda}$.  
  We have $\Lex_\Omega = \Lex'_\Omega$ if and only if 
  the chain $\Omega$ is built by adding one element at a time,
  i.e.\ $J_{\lambda+1} = J_\lambda \sqcup \{j_{\lambda+1}\}$ 
  for all $\lambda < \ell$.
\end{proposition}

\begin{proof}
  We will first establish the existence of a suitable chain $\Omega$.
  Given $C \colon \PPO_A^I \to \PPO_A$, we set $J_0 := \emptyset$ 
  and inductively define $J_{\lambda+1} := \Delta{J_\lambda}$.  
  Since $I$ is finite, this sequence will stabilize with some 
  $J_\ell = \Delta{J_\ell}$.  We thus obtain an ascending chain 
  $\Omega = \{ J_1 \subset \dots \subset J_\ell \}$.

  \begin{itemize}
  \item
    In order to show $\Lex'_\Omega \subseteq C$ consider 
    $J = J_\lambda$ and $K = J_{\lambda+1} \minus J_\lambda$.
    We have $K = \delta{J}$ by construction.  As a consequence,
    if $a \ndiff_J b$ and $a \spref_K b$ then $a \wpref b$.
    By monotonicity, the same holds true if $a \wpref_J b$ and $a \spref_K b$.
    This proves the inclusion $\Lex'_\Omega \subseteq C$.
  \item
    In order to show $C \subseteq \Lex_\Omega$ suppose that $a \wpref b$.
    The inclusion is obvious if $a \ndiff_J b$ for all $J \in \Omega$.
    Otherwise there exists $J = J_\lambda$ and $K = J_{\lambda+1} \minus J_\lambda$ 
    such that $a \ndiff_J b$ but not $a \ndiff_K b$.  Since $K$ is the minimal decisive 
    set relative to $J$, we have $a \wpref_K b$ by Proposition \ref{prop:VetoPower}.
    This proves the inclusion $C \subseteq \Lex_\Omega$.
  \end{itemize}
  
  Conversely, suppose that $\Omega = \{ J_1 \subset \dots \subset J_\ell \}$
  is a chain such that $\Lex'_\Omega \subseteq C \subseteq \Lex_\Omega$.
  The inclusion $\Lex'_\Omega \subseteq C$ implies that $J_1$ is decisive 
  for $C$, and $C \subseteq \Lex_\Omega$ implies that $J_1$ is minimal
  (see Proposition \ref{prop:VetoPower}).  This argument can be iterated: 
  $J_{\lambda+1} \minus J_\lambda$ is decisive relative to $J_\lambda$ 
  because $\Lex'_\Omega \subseteq C$, and minimal because $C \subseteq \Lex_\Omega$.
  We conclude that $J_1 = \Delta\emptyset$ and inductively
  $J_{\lambda+1} = \Delta{J_\lambda}$ for all $\lambda$.
\end{proof}

\subsection{Coalition structure of arrovian juntas} \label{sub:Classification}

The notion of relatively decisive subsets has allowed us 
to analyze examples such as $\Lex_\Omega$ and $\Lex'_\Omega$.  
As a final step, we will now use it to classify all arrovian 
voting systems for a finite society $I$:

\begin{theorem} \label{thm:PrincipalClassification}
  Given an arrovian voting system $C \colon \PPO_A^I \to \PPO_A$
  on three or more alternatives, there exists a map $\delta \colon 2^I \to 2^I$
  such that $\delta{N}$ is the minimal decisive set relative to $N$.
  Moreover, $\Delta \colon 2^I \to 2^I$, $\Delta{N} = \delta{N} \sqcup N$,
  is the unique map enjoying the following properties:
  \begin{enumerate}
  \item
    $N \subseteq \Delta{N}$ and $\Delta{M} \subseteq \Delta{N}$ 
    for all $M \subseteq N \subseteq I$.
    \hfill (Monotonicity)
  \item
    If $\Delta{N} \cap M = \emptyset$ for some $N,M \subseteq I$,
    then $\Delta( N \sqcup M ) = \Delta{N} \sqcup M$.
    \hfill (Minimality)
  \item
    We have $a \wpref b$ if and only if $a \ndiff_N b$ 
    and $a \spref_{\delta{N}} b$ for some $N \subseteq I$.
    \hfill (Decisiveness)
  \end{enumerate}
  
  Conversely, given an arbitrary map $\Delta \colon 2^I \to 2^I$ 
  satisfying (1), we set $\delta{N} = \Delta{N} \minus N$ and define 
  a map $C_\Delta \colon \PPO_A^I \to \PPO_A$ by the rule (3).
  The resulting voting system $C_\Delta$ satisfies 
  the arrovian axioms of unanimity, neutrality, monotonicity, 
  and independence of irrelevant alternatives.  Condition (2) 
  ensures that $\delta{N}$ is the minimal decisive set relative to $N$.
\end{theorem}

Before proving this theorem, let us clarify some details and reformulations:

\begin{lemma} \label{lem:WeakMinimality}
  If (1) holds, then (2) becomes equivalent to the following weaker condition:
  \begin{enumerate}
  \item[(2')]
    If $\Delta{N} \cap M = \emptyset$ for some $N,M \subseteq I$,
    then $\Delta( N \sqcup M ) \subseteq \Delta{N} \sqcup M$.
  \end{enumerate}
\end{lemma}

\begin{proof}
  Obviously (2) implies (2').  Given (1),
  $M \subseteq N \sqcup M$ implies $\Delta{M} \subseteq \Delta(N \sqcup M)$,
  and $N \subseteq N \sqcup M$ implies $\Delta{N} \subseteq \Delta(N \sqcup M)$.
  We obtain $\Delta{N} \cup \Delta{M} \subseteq \Delta(N \sqcup M)$,
  and in particular $\Delta{N} \sqcup M \subseteq \Delta(N \sqcup M)$.
  This proves that (1) and (2') imply (2).
\end{proof}

\begin{lemma} 
  If (1) holds, then (3) becomes equivalent to the following condition:
  \begin{enumerate}
  \item[(3')]
    We have $a \wpref b$ if and only if $a \wpref_N b$ 
    and $a \spref_{\delta{N}} b$ for some $N \subseteq I$.
  \end{enumerate}
\end{lemma}

\begin{proof}
  Obviously, if $a \ndiff_N b$ and $a \spref_{\delta{N}} b$,
  then $a \wpref_N b$ and $a \spref_{\delta{N}} b$.
  Conversely, if $a \wpref_N b$ and $a \spref_{\delta{N}} b$,
  then $a \ndiff_{N'} b$ and $a \spref_{K} b$
  for $N' = \{ i \in N \mid a \ndiff_i b \}$
  and $K = \delta{N} \sqcup (N \minus N')$.
  By (1) we have $\delta{N'} \subseteq \delta{N} \subseteq K$, 
  hence $a \ndiff_{N'} b$ and $a \spref_{\delta{N'}} b$ as desired.
\end{proof}

\begin{lemma} \label{lem:WeakDecisiveness}
  If (2) holds, then (3) becomes equivalent to the following condition:
  \begin{enumerate}
  \item[(3'')]
    We have $a \wpref b$ if and only if $a \spref_{\delta{N'}} b$
    for $N' = \{ i \in I \mid a \ndiff_i b \}$.
  \end{enumerate}
\end{lemma}

\begin{proof}
  We only have to show that $a \ndiff_N b$ and 
  $a \spref_{\delta{N}}$ for some $N \subseteq I$ imply 
  $a \spref_{\delta{N'}}$ for the possibly larger set $N'$.
  Obviously $M = N' \minus N$ is disjoint from $\delta{N}$.  
  Condition (2) now ensures that $\delta{N'} = \delta{N}$, 
  hence $a \spref_{\delta{N'}} b$ as desired.
\end{proof}

\begin{lemma} \label{lem:Decisiveness}
  If $a \wpref b$, then $K = \{ i \in I \mid a \spref_i b \}$
  is decisive relative to $N = \{ i \in I \mid a \ndiff_i b \}$.
\end{lemma}

\begin{proof}
  This is a consequence of Lemma \ref{lem:DecisiveSupport}:
  restricting to the electorate $I \minus N$ we have
  $K = \{ i \in I \minus N \mid a \wpref_i b \}$.
\end{proof}

\begin{proof}[Proof of Theorem \ref{thm:PrincipalClassification}]
  We already know from Lemma \ref{lem:NeutralityMonotonicity} that every arrovian 
  voting system $C \colon \PPO_A^I \to \PPO_A$ satisfies neutrality and monotonicity.  
  If $I$ is finite then we can construct the map $\delta \colon 2^I \to 2^I$, 
  $N \mapsto \delta{N}$, by applying Proposition \ref{prop:RestrictionOfElectorate}.
  Condition (3) holds by construction of $\delta$
  and the preceding Lemma \ref{lem:Decisiveness}.
  Turning to Condition (1), we first remark that 
  $N \subseteq \Delta{N}$ by definition of $\Delta$.
  Consider $M \subseteq N \subseteq I$: we know that
  $a \ndiff_N b$ and $a \spref_{\delta{N}} b$ imply $a \wpref b$.
  For $K = \Delta{N} \minus M$, monotonicity of $C$ ensures
  that $a \ndiff_M b$ and $a \spref_K b$ imply $a \wpref b$.
  This means that $\delta{M} \subseteq K$, hence 
  $\Delta{M} \subseteq \Delta{N}$, as claimed.
  For (2), we know that $a \ndiff_N b$ and $a \spref_{\delta{N}} b$
  imply $a \wpref b$, independently of all other preferences.  
  If $\Delta{N} \cap M = \emptyset$, then we can additionally
  assume $a \ndiff_M b$ and still conclude $a \wpref b$.  
  Thus $\delta(N \sqcup M) \subseteq \delta{N}$, or equivalently, 
  $\Delta(N \sqcup M) \subseteq \Delta{N} \sqcup M$.
  For the inverse inclusion see Lemma \ref{lem:WeakMinimality}.

  In order to prove uniqueness we have to show, for every 
  map $\Delta \colon 2^I \to 2^I$ satisfying (1), (2), (3),
  that $\delta{N} = \Delta{N} \minus N$ is the smallest
  decisive set relative to $N$.  To begin with, Condition (3) 
  states that $\delta{N}$ is decisive relative to $N$.  
  Conversely, consider an arbitrary set $K \subseteq I \minus N$ 
  that is decisive relative $N$.  By definition, 
  $a \ndiff_{N} b$ and $a \spref_{K} b$ imply $a \wpref b$, 
  even if $a \sferp_J b$ on $J = I \minus ( N \sqcup K )$.
  According to Condition (3) there exists $N' \subseteq I$
  such that $a \ndiff_{N'} b$ and $a \spref_{\delta{N'}} b$.
  We necessarily have $N' \subseteq N$ and $\delta{N'} \subseteq K$.  
  Since $M = N \minus N'$ is disjoint from $\delta{N'}$, 
  Condition (2) implies that $\delta{N} = \delta{N'}$.
  We conclude that $\delta{N} \subseteq K$, which means that 
  $\delta{K}$ is indeed the smallest decisive set relative to $N$.

  Conversely, given $\Delta \colon 2^I \to 2^I$ satisfying (1), 
  we have to verify that the associated voting system 
  $C_\Delta \colon \PPO_A^I \to \PPO_A$ is well-defined:
  the outcome is obviously a reflexive relation, 
  and the only delicate point is transitivity.
  Given $a \wpref b$ and $b \wpref c$ we want to show $a \wpref c$.
  We know that $a \ndiff_{N_1} b$ and $a \spref_{\delta{N_1}} b$,
  as well as $b \ndiff_{N_2} c$ and $b \spref_{\delta{N_2}} c$,
  for some $N_1, N_2 \subseteq I$.  We thus have $a \ndiff_N c$
  on the intersection $N = N_1 \cap N_2$, and $a \spref_K c$ 
  on $K = (\Delta{N_1} \cap \Delta{N_2}) \minus N$.
  Condition (1) ensures that $\Delta{N} \subseteq \Delta{N_1}$ 
  and $\Delta{N} \subseteq \Delta{N_2}$, thus 
  $\Delta{N} \subseteq \Delta{N_1} \cap \Delta{N_2}$.
  This proves that $\delta{N} \subseteq K$, whence $a \spref_{\delta{N}} c$.
  The defining rule (3) now says that $a \wpref c$, as desired.

  Having proved that the map $C_\Delta \colon \PPO_A^I \to \PPO_A$ 
  is well-defined, we conclude that it satisfies the arrovian axioms:
  unanimity, neutrality, and independence of irrelevant alternatives 
  are clear by (3), while monotonicity follows from (1).
  Finally, it remains to show for $C_\Delta$ that 
  $\delta{N} = \Delta{N} \minus N$ is indeed 
  the smallest decisive set relative to $N$:
  given conditions (1), (2), and (3), this 
  follows from the uniqueness proved above.
\end{proof}

\begin{remark}
  The theorem asserts that every arrovian voting system 
  $C \colon \PPO_A^I \to \PPO_A$ on three or more alternatives 
  is characterized by the associated map $\Delta \colon 2^I \to 2^I$. 
  In particular we can extract the usual global information:
  $K = \Delta\emptyset$ is the smallest decisive subset, 
  whereas the union $J = \bigcup_\lambda \Delta^\lambda K$ 
  is the smallest strongly decisive subset.
  More precisely, we can recover the chain
  $\Omega = \{ J_1 \subset \dots \subset J_\ell \}$
  defined by $J_1 = K$ and $J_{\lambda+1} := \Delta{J_\lambda}$. 
  According to Proposition \ref{prop:LexicographicInclusion}
  we have $\Lex'_\Omega \subseteq C \subseteq \Lex_\Omega$
  but in general the inclusion can be strict.  
\end{remark}

\begin{example}
  The smallest voting systems that are not 
  of lexicographic type occur for three voters.
  Example \ref{exm:NonLexicographic} illustrates that the 
  classification is not exhausted by lexicographic voting systems alone. 
  Here we find the chain $\Omega = \{ \{1,2\} \subset \{1,2,3\} \}$.
  We have $C \ne \Lex_\Omega$ because $\Delta_C\{1\} = \{1,2,3\}$ but $\Delta_\Omega\{1\} = \{1,2\}$.
  Likewise $C \ne \Lex'_\Omega$ because $\Delta_C\{2\} = \{1,2\}$ but $\Delta'_\Omega\{2\} = \{1,2,3\}$.
  We conclude that $\Lex'_\Omega \subsetneq C \subsetneq \Lex_\Omega$.
\end{example}

\begin{corollary}
  We have $C \subseteq C'$ if and only if $\Delta \supseteq \Delta'$.
  More explicitely, $C(P_1,\dots,P_n) \subseteq C'(P_1,\dots,P_n)$
  for all profiles $(P_1,\dots,P_n) \in \PPO_A^I$ if and only if
  $\Delta{N} \supseteq \Delta'{N}$ for all $N \subseteq I$.
\end{corollary}

\begin{proof}
  We first recall that $\Delta{N} = \delta{N} \sqcup N$ and 
  $\Delta'{N} = \delta'{N} \sqcup N$, which means that
  the conditions $\Delta{N} \supseteq \Delta'{N}$ and 
  $\delta{N} \supseteq \delta'{N}$ are equivalent.
  
  $(\Rightarrow)$
  Given $C \subseteq C'$ we want to show that every subset 
  $K \subseteq I$ that is decisive for $C$ relative to $N$
  is also decisive for $C'$ relative to $N$.  Consider $a \ndiff_N b$ and 
  $a \spref_K b$ and $a \sferp_L b$ on the complement $L = I \minus (N \sqcup K)$.
  Since $K$ is decisive relative to $N$ we know that $a \wpref b$,
  and from $C \subseteq C'$ we deduce $a \wpref' b$.  
  We conclude that $K$ is decisive for $C'$ relative to $N$.

  $(\Leftarrow)$
  We suppose that every subset $K \subseteq I$ that is decisive 
  for $C$ relative to $N$ is also decisive for $C'$ relative to $N$.
  If $a \wpref b$ then $K = \{ i \mid a \spref_i b \}$ 
  is decisive for $C$ relative to $N = \{ i \mid a \ndiff_i b \}$.
  By hypothesis, $K$ is also decisive for $C'$ relative to $N$.
  We conclude that $a \wpref' b$, as claimed.
\end{proof}

\subsection{Back to linear orderings} \label{sub:Applications}

The preceding proofs show a little more than stated in 
the Classification Theorem \ref{thm:PrincipalClassification}:
one can weaken the requirements by demanding that the voting system $C$ 
be only defined on the subset $\LPO_A \subset \PPO_A$ 
of linear orderings (i.e.\ complete preorders as opposed to partial preorders).
We thus obtain the following result:

\begin{corollary} \label{cor:LinearClassification}
  For every arrovian voting system $C \colon \LPO_A^I \to \PPO_A$
  on three or more alternatives there exists a unique map 
  $\Delta \colon 2^I \to 2^I$ satisfying the conditions 
  of Theorem \ref{thm:PrincipalClassification}
  such that $C = C_\Delta|\LPO_A^I$.
\end{corollary}

\begin{proof}
  As explained in the proof of Proposition \ref{prop:RestrictionToLinearOrderings},
  the arguments developed in \textsection\ref{sec:DecisiveSubsets} for partial 
  orderings also apply to linear orderings, and all constructions 
  can be carried out within $\LPO_A^I \subset \PPO_A^I$.
  We can then consider relatively decisive subsets 
  (Definition \ref{def:RelativelyDecisiveSubset}), and 
  the argument of restricting the electorate works as before 
  (Proposition \ref {prop:RestrictionOfElectorate}).
  In particular, for each $N \subseteq I$ there exists a unique minimal 
  decisive subset relative to $N$, denoted $\delta{N} \subseteq I \minus N$.
  The proof of Theorem \ref{thm:PrincipalClassification} then
  translates verbatim to a proof of Corollary \ref{cor:LinearClassification}.
\end{proof}

One could also strengthen the requirements and demand 
that $C$ take values in the set of linear orderings.
Here is a typical example:

\begin{example} \label{exm:StrictLexicographical}
  We consider a family $\pi=(k_1,\dots,k_\ell)$ of distinct elements of $I$ 
  and define $\Lex_{\pi} \colon \LPO_A^I \to \LPO_A$ as follows:
  We set $a \ndiff b$ if and only if $a \ndiff_k b$ for each $k=k_1,\dots,k_\ell$.
  Otherwise, let $k$ be the first element in the family 
  $k_1,\dots,k_\ell$ with $a \not\ndiff_k b$,
  and set either $a \spref b$ if $a \spref_k b$,
  or $a \sferp b$ if $a \sferp_k b$.
  This is the \emph{lexicographic} voting rule
  associated to the family $\pi=(k_1,\dots,k_\ell)$.
  As indicated, $k_\lambda$ has priority over $k_{\lambda+1}$, and only 
  in the case of indifference is the decision passed down in the hierarchy.
  Notice that $K = \{k_1\}$ is the minimal decisive set, whereas 
  $J = \{k_1,\dots,k_\ell\}$ is the minimal strongly decisive set,
  and $\pi$ is a permutation of the set $J$.
\end{example}

Although the change of domain from $\PPO_A^I$ to $\LPO_A^I$ 
has turned out to be insignificant, changing the range 
from $\PPO_A$ to $\LPO_A$ alters the classification result 
dramatically.  We are led to the following refinement of Arrow's 
dictator theorem, as stated in Theorem \ref{Thm:ArrowRefined}:

\begin{corollary}[Refined version of Arrow's dictator theorem]
  \label{cor:DictatorTheorem}
  Suppose that the set $A$ contains at least three alternatives.
  Then every map $C \colon \LPO_A^I \to \LPO_A$ that satisfies 
  unanimity and IIA is a lexicographic voting rule of the form 
  $C = \Lex_{\pi}$ for a unique family $\pi=(k_1,\dots,k_\ell)$.
\end{corollary}

\begin{proof}
  Via the inclusion $\LPO_A \subset \PPO_A$,
  we can apply the previous classification theorem
  to the map $C \colon \LPO_A^I \to \PPO_A$.
  For each $N \subseteq I$, the minimal relatively decisive 
  subset $\delta{N}$ contains at most one element.
  If it contained two different elements $i,j \in K$, 
  then $a \ndiff_N b$ together with $a \spref_i b$ and $a \sferp_j b$ 
  would imply $a \ncomp b$, and $C$ would not map to $\LPO_A$.
  (Recall that $i$ and $j$ have veto power,
  see Proposition \ref{prop:VetoPower}.)

  Let $K = \delta\emptyset$ be the minimal decisive subset.
  If $K = \emptyset$ then $C$ is the trivial voting system 
  and we conclude that $\ell=0$.  Otherwise $K$ contains exactly 
  one element, $K = \{k_1\}$, with $k_1$ being the dictator. 
  In this case, $\delta\{k_1\}$ is either empty or contains 
  exactly one element, $\delta\{k_1\} = \{k_2\}$.
  In the latter case we find $\delta\{k_1,k_2\} = \emptyset$
  or $\delta\{k_1,k_2\} = \{k_3\}$.  Iterating this argument, 
  we obtain a family $\pi=(k_1,\dots,k_\ell)$ of distinct elements of $I$.
  We conclude that $C = \Lex_{\pi}$ according to Proposition \ref{prop:LexicographicInclusion}.
\end{proof}

\begin{remark}
  If we insist that the voting system be strongly unanimous,
  then the only solutions are of the form $\Lex_{\pi}$ where 
  $\pi = (k_1,\dots,k_n)$ is a permutation of the entire set $I$.
  There are thus precisely $n!$ voting systems $C \colon \LPO_A^I \to \LPO_A$
  that are strongly unanimous and independent of irrelevant alternatives.
  Even though dictatorial, they have the advantage to extract a maximum 
  of information within the arrovian setting.  The only point of unfairness,
  of course, is the arbitrary choice of the permutation $\pi$, 
  that is, the order of individuals in the hierarchy.
\end{remark}


\section{Infinite societies} \label{sec:InfiniteSocieties}

We conclude this article by adapting our arguments 
to the case where the set $I$ of voters is infinite.
The results are more involved but nevertheless illuminating 
by placing the finite case in a wider perspective.
As before, $\PPO_A^I$ denotes the set of all
preference profiles $(P_i)_{i \in I}$, or equivalently, 
of all maps $I \to \PPO_A$, $i \mapsto P_i$.

\subsection{Lexicographic voting rules} \label{sub:InfiniteLexicographic}

For infinite societies, lexicographic voting rules can be defined 
as in Examples \ref{exm:Lexicographic} and \ref{exm:StrongLexicographic}.
Here we consider a chain $\Omega \subseteq 2^I$, i.e.\ a collection 
of subsets of $I$ that is linearly ordered by inclusion.
The only subtlety is that $\Omega$ must be \emph{well-ordered}, 
that is, every non-empty subset $\Omega' \subseteq \Omega$ has a minimal element,
i.e.\ the intersection $\bigcap_{J \in \Omega'} J$ is again an element of $\Omega'$.
(This is automatically satisfied if $\Omega$ is finite.)

\begin{proposition}[Lexicographic voting rule] \label{prop:WeakLexicographic}
  If $\Omega \subseteq 2^I$ is a well-ordered chain, then
  the voting systems $\Lex_\Omega, \Lex'_\Omega \colon \PPO_A^I \to \PPO_A$
  are well-defined and satisfy the axioms of unanimity, neutrality, 
  monotonicity, and independence of irrelevant alternatives.
  \qed
\end{proposition}

As a special case, consider a subset $J \subseteq I$ equipped 
with some well-ordering $\le$.  Then the initial segments
$J_j = \{ i \in J \mid i \le j \}$ form a well-ordered 
chain $\Omega = \{ J_j \mid j \in J \}$.  The associated voting 
rule $\Lex_{(J,\le)} = \Lex_\Omega = \Lex'_\Omega$ works as follows:
we have $a \ndiff b$ if and only if $a \ndiff_J b$.  Otherwise
let $j \in J$ be the smallest element such that $a \not\ndiff_j b$;
we then set $a \wpref b$ if and only if $a \wpref_j b$, and 
symmetrically $a \wferp b$ if and only if $a \wferp_j b$.
As a concrete example, consider 
$\Lex_{(\N,\le)} \colon \PPO_A^\N \to \PPO_A$ 
where $J = \N$ is ordered as usual, thus
$\Omega = \{ \{1\} \subset \{1,2\} \subset \{1,2,3\} \subset \dots \}$.

\begin{proposition}[Strict lexicographic voting rule] \label{prop:StrictLexicographic}
  Suppose that $J \subseteq I$ is equipped with some well-ordering $\le$.
  Then the voting system $\Lex_{(J,\le)} \colon \PPO_A^I \to \PPO_A$ 
  maps $\LPO_A^I$ to $\LPO_A$. 
  Moreover, for each $N \subseteq I$ there exists 
  a unique minimal decisive subset $\delta{N} \subseteq I \minus N$: 
  we have $\delta{N} = \emptyset$ if and only if $J \subseteq N$,
  and otherwise $\delta{N} = \{ \min(J \minus N) \}$.
  \qed
\end{proposition}

\subsection{Principal voting systems} \label{sub:PrincipalVotingSystems}

The uniqueness of a minimal decisive subset $K \subseteq I$ 
remains valid even if the set $I$ of voters is infinite:  
if two subsets $K_1,K_2 \subseteq I$ are decisive 
and minimal, then their intersection is decisive, 
hence $K_1 = K_1 \cap K_2 = K_2$ by minimality.

Our hypothesis that the set $I$ be finite is crucial, however, 
for establishing the \emph{existence} of a minimal decisive set:

\begin{example} \label{exm:CofiniteVotingRule}
  Consider an infinite set $I$ 
  and define $C \colon \PPO_A^I \to \PPO_A$ by the rule $a \wpref b$ 
  if and only if $a \wpref_i b$ for all but finitely many voters  $i \in I$.
  Here every finite set $J \subseteq I$ is negligible since it has 
  no influence on the outcome.  Conversely, a subset $J \subseteq I$
  is decisive if and only if it has finite complement $I \minus J$.
  There is, however, no minimal such set.
\end{example}

\begin{definition} \label{def:PrincipalVotingSystems}
  An arrovian voting system $C \colon \PPO_A^I \to \PPO_A$ 
  is called \emph{principal} if for each $N \subseteq I$ there exists 
  a minimal subset $M \subseteq I \minus N$ that is decisive relative to $N$.
  In this case $M$ is uniquely determined by $N$ 
  and will be denoted by $\delta{N}$ as before.
\end{definition}

\begin{remark} \label{rem:PrincipalVotingSystems}
  Every arrovian voting system on a finite set of voters $I$ is principal.
  For an infinite society, $C \colon \PPO_A^I \to \PPO_A$ can be principal, 
  for example $\Lex_{(J,\le)}$, or non-principal, as in 
  the preceding Example \ref{exm:CofiniteVotingRule}.
  For principal voting systems the Classification Theorem
  \ref{thm:PrincipalClassification} holds verbatim: given the existence 
  of $\delta{N}$, the finiteness of $I$ is not used in the proof.
  Likewise we have the variant for linear orderings
  formulated in Corollary \ref{cor:LinearClassification}.
\end{remark}

\begin{definition}
  Every well-ordered chain $\Omega$ of subsets of $I$
  can be uniquely indexed by ordinal numbers $\lambda$ 
  such that $\Omega = \{ J_\lambda \mid \lambda \le \ell \}$
  and $J_\kappa \subset J_\lambda$ for $\kappa < \lambda$.
  We say that $\Omega$ is \emph{continuously} well-ordered
  if  $J_\lambda = \bigcup_{\kappa < \lambda} J_\kappa$ 
  for every limit ordinal $\lambda$.
\end{definition}

\begin{proposition} \label{prop:InfiniteLexicographicInclusion}
  For every principal arrovian voting system $C \colon \PPO_A^I \to \PPO_A$ 
  there exists a unique continuously well-ordered chain $\Omega$ 
  such that $\Lex'_\Omega \subseteq C \subseteq \Lex_\Omega$, namely the chain 
  inductively defined by $J_0 := \emptyset$, and $J_{\lambda+1} := \Delta{J_\lambda}$. 
  We have $\Lex_\Omega = \Lex'_\Omega$ if and only if the chain $\Omega$ is built 
  by adding one element at a time, i.e.\ $J_{\lambda+1} = J_\lambda \sqcup \{j_{\lambda+1}\}$ 
  for all $\lambda < \ell$.
\end{proposition}

\begin{proof}
  By hypothesis, $C$ is principal, so for each $N \subseteq I$ 
  we can consider the minimal relatively decisive set $\delta{N}$.
  We set $\Delta{N} = \delta{N} \sqcup N$ and define $J_0 := \emptyset$,
  and inductively $J_{\lambda+1} := \Delta{J_\lambda}$.
  In the case of an infinite set $I$, we proceed by transfinite 
  induction, setting $J_\lambda := \bigcup_{\kappa < \lambda} J_\kappa$ 
  for every limit ordinal $\lambda$.  Since $I$ is a set, this process must
  stop with $\Delta{J_\ell} = J_\ell$ for some ordinal $\ell$, and we obtain 
  a continuously well-ordered chain $\Omega = \{ J_\lambda \mid \lambda \le \ell \}$.
  The double inclusion $\Lex'_\Omega \subseteq C \subseteq \Lex_\Omega$
  follows as in the proof of Proposition \ref{prop:LexicographicInclusion},
  and the uniqueness argument generalizes verbatim 
  to continuously well-ordered chains.
\end{proof}

We are now in position to prove the converse of 
Proposition \ref{prop:StrictLexicographic} and 
thus characterize strict lexicographic voting rules.
We obtain the following refined version of Arrow's dictator 
theorem, which comprises the finite and the infinite case:

\begin{theorem} \label{thm:InfiniteDictators}
  Suppose that $C \colon \PPO_A^I \to \PPO_A$ is a principal
  arrovian voting system that maps $\LPO_A^I$ to $\LPO_A$.
  Then we have $C = \Lex_{(J,\le)}$ for some subset 
  $J \subseteq I$ equipped with a well-ordering $\le$, 
  and the pair $(J,\le)$ is uniquely determined by $C$.

  Conversely, every pair $(J,\le)$ can be realized in this way,
  thus establishing a bijection between the said voting systems
  and well-ordered sets $(J,\le)$ with $J \subseteq I$.

  Moreover, strong unanimity is equivalent to $J = I$.
  Principal arrovian voting system $C \colon \PPO_A^I \to \PPO_A$
  that map $\LPO_A^I$ to $\LPO_A$ and satisfy strong unanimity
  are thus in bijective correspondence with well-orderings 
  of the electorate $I$.
\end{theorem}

\begin{proof} 
  By Proposition \ref{prop:InfiniteLexicographicInclusion}, 
  we have $\Lex'_\Omega \subseteq C \subseteq \Lex_\Omega$ for a unique continuously 
  well-ordered chain $\Omega$, but in general this inclusion may be strict.
  We will now exploit the hypothesis that $C(\LPO_A^I) \subseteq \LPO_A$.
  As shown in the proof of Corollary \ref{cor:DictatorTheorem}, 
  each $\delta{N}$ is either empty or consists of a single individual.
  The set $J = J_\ell$ becomes well-ordered via the bijection 
  $\lambda \mapsto j_\lambda$ defined by the condition
  $J_{\lambda+1} \minus J_\lambda = \{ j_{\lambda+1} \}$.
  According to Proposition \ref{prop:InfiniteLexicographicInclusion},
  we conclude that $C = \Lex_{(J,\le)}$ and the pair $(J,\le)$ is 
  uniquely determined by $C$.
\end{proof}

\subsection{The filter of decisive subsets} \label{sub:Filters}

For non-principal voting systems we will now explain how to generalize 
the Classification Theorem \ref{thm:PrincipalClassification}.
As P.C.\,Fishburn \cite{Fishburn:1970} pointed out,
the family $\F$ of decisive sets forms a \emph{filter} 
in the following sense: 

\begin{definition}
  A \emph{filter} on a set $I$ is a collection 
  of subsets $\F \subseteq 2^I$ such that
  \begin{enumerate}
  \item[(F1)] 
    If $K \subseteq J \subseteq I$, then $K \in \F$ implies $J \in \F$.
  \item[(F2)] 
    We have $I \in \F$, that is, $\F$ is non-empty.
  \item[(F3)] 
    If $K,J \in \F$, then $K \cap J \in \F$.
  \end{enumerate}

  Given $K \subseteq I$ the collection
  $(K) = \{ J \subseteq I \mid K \subseteq J \}$ 
  is the \emph{principal} filter generated by $K$.
  One has $\emptyset \in \F$ if and only if $\F$ is 
  the \emph{trivial} filter, i.e.\ $\F = (\emptyset) = 2^I$.
  A filter $\F$ is called \emph{proper} if $\emptyset \notin \F$.
  An \emph{ultrafilter} on $I$ is a maximal proper filter.
\end{definition}

\begin{remark}
  Following Bourbaki \cite[\textsection I.6.1]{Bourbaki:1971},
  most authors demand $\emptyset \notin \F$ as a fourth filter axiom.
  We will not do so here because trivial filters naturally occur
  in the sequel, as the decisive sets of trivial voting systems.
  Moreover, a filter on a set $I$ is analogous to an ideal in a ring $R$, 
  and this definition usually includes the ring itself as the trivial ideal.%
  \footnote{ An ideal $S$ in a ring $(R,+,\cdot)$ 
    is a subset $S \subseteq R$ such that 
    (i) $r \cdot s \in S$ for all $r \in R$ and $s \in S$,
    (ii) $0 \in S$, that is, $S$ is non-empty, and 
    (iii) $s + t \in S$ for all $s,t \in S$.
    If we replace the ring $(R,+,\cdot)$ by the boolean algebra $(2^I,\cap,\cup)$,
    then the conditions (i),(ii),(iii) translate to the filter axioms (F1),(F2),(F3).  
    Notice that $0 \in R$ is the neutral element with respect to addition $+$,
    while $I \in 2^I$ is the neutral element with respect to intersection $\cap$. }
  If this is to be excluded, one should speak of \emph{proper} ideals.  
  We will conform our notation to this algebraic analogy, and thus 
  speak of \emph{proper} filters if we wish to exclude the trivial case.
\end{remark}

\begin{remark}
  If a filter $\F$ contains a minimal element $K$, then $K$ 
  is unique and $\F = (K)$ is the principal filter generated by $K$.
  If $I$ is infinite then a filter $\F$ does not necessarily have 
  a minimal element: consider the filter $\F$ of all cofinite subsets 
  $K \subseteq I$, that is, subsets $K$ with finite complement $I \minus K$.
  (This is called the Fr\'echet filter on $I$.) 
\end{remark}

\begin{remark}
  Given an element $i \in I$, the principal filter $(\{i\})$ is an ultrafilter.
  If $I$ is finite, then every ultrafilter $\F$ is of the form $\F = (\{i\})$.
  In general, the axiom of choice guarantees that 
  every proper filter is contained in some ultrafilter.  
  A filter $\F$ is an ultrafilter if and only if for every subset 
  $K \subseteq I$ one has either $K \in \F$ or $I \minus K \in \F$.
  \cite[\textsection I.6.4]{Bourbaki:1971}
\end{remark}

Given a filter $\F$ we define $a \wpref_\F b$ 
to signify $a \wpref_J b$ for some $J \in \F$.
In the case of a principal filter $\F = (K)$
we have that $a \wpref_{(K)} b$ is equivalent 
to $a \wpref_K b$, as defined previously.
Analogously we define $a \spref_\F b$ 
to signify $a \spref_J b$ for some $J \in \F$, etc.

The following result was first published by A.P.\,Kirman 
and D.\,Sondermann \cite{KirmanSondermann:1972}, and implicitly
in P.C.\,Fishburn's previous article \cite{Fishburn:1970}.
Essentially, it had already been discovered in 1952 
by G.\,Guilbaud \cite{Guilbaud:1952}.  In order to make
our presentation self-contained, we will state and prove 
the essential observation needed for our classification:

\begin{proposition}[Filter of decisive subsets] \label{prop:InfiniteFilter}
  For every arrovian voting system $C \colon \PPO_A^I \to \PPO_A$ 
  on three or more alternatives, the family $\F$ of decisive sets is 
  a filter on $I$.  The same holds true for the subfamily $\F' \subseteq \F$ 
  of strongly decisive sets.  Moreover, if $C$ maps $\LPO_A^I$ to $\LPO_A$, 
  then $\F$ is either trivial or an ultrafilter on $I$.

  Conversely, every filter $\F$ on $I$ defines an arrovian voting system
  $C_\F \colon \PPO_A^I \to \PPO_A$ by setting $a \wpref b$ if and only if 
  $a \wpref_\F b$, and in this case $\F$ is the filter of decisive sets
  and also the filter of strongly decisive sets.
  If $\F = (K)$ is principal, then $C_{(K)} = C_K$ as defined previously.
  If $\F$ is either trivial or an ultrafilter on $I$, 
  then $C_\F$ maps $\LPO_A^I$ to $\LPO_A$.
\end{proposition}

\begin{proof}
  For every arrovian voting system $C$, the family $\F$ resp.\ $\F'$
  is a filter: it satisfies (F1) by monotonicity, (F2) by unanimity, 
  and (F3) by Proposition \ref{prop:IntersectionProperty}.

  Suppose, moreover, that $C(\LPO_A^I) \subseteq \LPO_A$.
  Consider a subset $J \subseteq I$ and its complement $K = I \minus J$.
  For a profile with $a \spref_J b$ and $a \sferp_K b$ we have two possible outcomes: 
  according to Remark \ref{rem:WorstCase}\,\eqref{item:WorstCaseDecisive},
  if $a \wpref b$ then $J$ is decisive; if $a \wferp b$ then $K$ is decisive.
  If we had both, then their intersection $J \cap K = \emptyset$ would be 
  decisive and $C$ would be the trivial voting system.  
  This being excluded, we conclude that $\F$ is an ultrafilter.

  Conversely, a filter $\F$ allows to define a map $C_\F \colon \PPO_A^I \to \PPO_A$.
  First of all we have to show that this is well-defined.
  Obviously, the outcome is a reflexive relation because $I \in \F$.
  The only delicate point is transitivity: given $a \wpref b$ 
  and $b \wpref c$, we know that $a \wpref_{K_1} b$ 
  and $b \wpref_{K_2} c$ for some $K_1, K_2 \in \F$. 
  We thus have $a \wpref_K c$ on the intersection $K = K_1 \cap K_2$.
  Condition (F3) ensures that $K \in \F$, and hence $a \wpref c$, as desired.

  Having proved that $C_\F \colon \PPO_A^I \to \PPO_A$ 
  is well-defined, we conclude that it satisfies 
  the arrovian axioms: neutrality, monotonicity and 
  independence of irrelevant alternatives are clear by construction.  
  Unanimity follows since $I \in \F$, as ensured by condition (F2).
  By definition, each set $J \in \F$ is decisive for the voting system $C_\F$.
  Conversely a subset $J \subseteq I$ is decisive if and only if 
  $K \subseteq J$ for some $K \in \F$.  Hence Condition (F1) 
  ensures that $\F$ is the family of decisive sets for $C_\F$.
 
  Suppose, moreover, that $\F$ is an ultrafilter.  Given two alternatives
  $a,b \in A$ consider the set $J = \{ i \in I \mid a \wpref_i b \}$ 
  and its complement $K = \{ i \in I \mid a \sferp_i b \}$.  
  Since $\F$ is an ultrafilter we have either $J \in \F$ or $K \in \F$.
  This shows that $a \wpref b$ or $a \wferp b$, 
  in other words, the outcome is a complete ordering.
\end{proof}

As in the finite case, it is easiest to classify voting systems
for which decisive and strongly decisive subsets co\"incide:

\begin{proposition}[Classification of juntas without internal structure] \label{prop:InfiniteSimpleJuntas}
  Given an arrovian voting system $C \colon \PPO_A^I \to \PPO_A$
  on three or more alternatives, let $\F$ be the filter of decisive sets, 
  and let $\F' \subseteq \F$ be the filter of strongly decisive sets.  
  We have the double inclusion $C_{\F'} \subseteq C \subseteq C_\F$,
  and equality $C = C_\F$ or $C = C_{\F'}$ holds if and only if $\F = \F'$, 
  that is, each decisive set is also strongly decisive.
\end{proposition}

\begin{proof}
  This is a variation of Proposition \ref{prop:SimpleJuntas}.
  The inclusion $C_{\F'} \subseteq C$ is clear
  by definition of strong decisiveness.
  In order to show $C \subseteq C_\F$ we appeal
  to Lemma \ref{lem:DecisiveSupport}:  given $a \wpref b$, 
  the supporting set $K = \{ i \in I \mid a \wpref_i b \}$ 
  is decisive, hence $a \wpref_K b$ with $K \in \F$.
\end{proof}

\subsection{Relatively decisive subsets} \label{sub:InfiniteClassification}

The preceding Proposition \ref{prop:InfiniteSimpleJuntas} characterizes 
voting systems in which decisive and strongly decisive subsets co\"incide.
In general they differ, as shown by lexicographic voting rules
(see \textsection\ref{sub:InfiniteLexicographic} above).
In order to classify all possibilities, we thus take up the detailed analysis
and consider the filter $\d{N}$ on $I \minus N$ of decisive subsets 
relative to $N$, following Proposition \ref{prop:RestrictionOfElectorate}.
If $I$ is finite then this is simply the principal filter $\d{N} = (\delta{N})$,
but in the infinite case $\d{N}$ may not be principal,
so that the language of filters is appropriate.
We can now reformulate the principal classification, 
Theorem \ref{thm:PrincipalClassification}, by replacing 
the set $\delta{N}$ with the filter $\d{N}$, which leads
to Theorem \ref{thm:MeasurableClassification} stated below.

We will end this \emph{tour de force} in set-theoretic 
abstraction by adding one final level of technicality.
For infinite sets $I$ it is sometimes inappropriate to consider 
arbitrary subsets $K \subseteq I$, that is, it may be necessary 
to work with some restricted family $\Sigma \subseteq 2^I$.
(See T.E.\,Armstrong \cite{Armstrong:1980,Armstrong:1985}.) 
Typically this occurs when $(I,\Sigma,\mu)$ is a measure space: 
quite often the measure $\mu \colon \Sigma \to \R_+$ is defined 
only on $\Sigma$ because it cannot be extended to the whole set $2^I$.
Consider for example $I = \R$ and $\mu \colon \Sigma \to \R_+$ 
the Lebesgue-measure defined on the family $\Sigma$ of 
Lebesgue-measurable sets.  Here the axiom of choice 
implies that $\Sigma \ne 2^\R$. 

All that has been said and done in this article 
generalizes in an obvious way to the measurable context. 
To be explicit, we demand $\Sigma$ to be an \emph{algebra} 
in the following sense, and that all subsets and filters 
respect this algebra:

\begin{definition}
  An \emph{algebra} on a set $I$ is a collection 
  of subsets $\Sigma \subseteq 2^I$ such that
  \begin{enumerate}
  \item We have $I \in \Sigma$, and $K \in \Sigma$ implies $I \minus K \in \Sigma$.
  \item If $K,J \in \Sigma$, then $K \cup J$ and $K \cap J$ are elements of $\Sigma$.
  \end{enumerate}
  The elements of $\Sigma$ are called \emph{measurable sets},
  and the pair $(I,\Sigma)$ is called a \emph{measurable space}.
  A \emph{filter} in an algebra $\Sigma$ is a subset $\F \subseteq \Sigma$ such that
  \begin{enumerate}
  \item[(F1)] 
    If $K,J \in \Sigma$, then $K \in \F$ and $K \subseteq J$ imply $J \in \F$.
  \item[(F2)] 
    We have $I \in \F$, that is, $\F$ is non-empty.
  \item[(F3)] 
    If $K,J \in \F$, then $K \cap J \in \F$.
  \end{enumerate}

  Given a measurable space $(I,\Sigma)$ and a set $K \in \Sigma$,
  the collection $(K) = \{ J \in \Sigma \mid K \subseteq J \}$ 
  is called the \emph{principal} filter generated by $K$.
  One has $\emptyset \in \F$ if and only if 
  $\F$ is the \emph{trivial} filter, i.e.\ $\F = (\emptyset) = \Sigma$.
  A filter $\F$ is called \emph{proper} if $\emptyset \notin \F$.
  An \emph{ultrafilter} in $\Sigma$ is a maximal proper filter in $\Sigma$.
\end{definition}

We define $\PPO_A^{\smash{I,\Sigma}}$ to be the family 
of measurable preference profiles $(P_i)_{i \in I}$, that is,
we demand the set $\{ i \in I \mid a \wpref_i b \}$ to be measurable 
for each pair of alternatives $a,b \in A$.  Since $\Sigma$ 
is an algebra, all relevant subsets of $I$ thus become measurable,
such as $\{ i \in I \mid a \ndiff_i b \}$,
or $\{ i \in I \mid a \spref_i b \}$, 
or $\{ i \in I \mid a \ncomp_i b \}$, etc.

A voting system for the society $(I,\Sigma)$ 
is a map $C \colon \PPO_A^{\smash{I,\Sigma}} \to \PPO_A$,
and the arrovian axioms can be formulated as before.
(Notice that $\Sigma = 2^I$ corresponds to 
a set $I$ without any measurability restrictions.)
If there are at least three alternatives, then unanimity 
and independence of irrelevant alternatives imply neutrality 
and monotonicity (Lemma \ref{lem:NeutralityMonotonicity}).
Moreover, the (strongly) decisive subsets $K \in \Sigma$ 
form a filter in the algebra $\Sigma$, and the preceding
Propositions \ref{prop:InfiniteFilter} and 
\ref{prop:InfiniteSimpleJuntas} still hold.

As in Proposition \ref{prop:RestrictionOfElectorate},
it is possible to restrict the electorate to $I \minus N$
for every measurable set $N \in \Sigma$, by passing from 
$\Sigma$ to $\Sigma^N = \{ K \minus N \mid K \in \Sigma \}$,
the restricted algebra on $I \minus N$.  As before, 
this trick allows to define $\d{N} \subseteq \Sigma^N$, 
the filter of decisive subsets relative to $N$.

\subsection{Voting systems for measurable societies} \label{sub:MeasurableSocieties}

Armed with the appropriate notation, 
the Classification Theorem 
now translates to arrovian voting systems 
$C \colon \PPO_A^{\smash{I,\Sigma}} \to \PPO_A$
for the measurable society $(I,\Sigma)$.
For ease of notation we define a map $\D \colon \Sigma \to \FF{\Sigma}$,
where $\FF{\Sigma}$ is the set of filters $\F \subseteq \Sigma$,
by $\D{N} = \d{N} + N = \{ J \sqcup N \mid J \in \d{N} \} \subseteq (N)$.
Since $\d{N} \subseteq \Sigma^N$, we can recover the initial data 
via $\d{N} = \D{N} - N = \{ J \minus N \mid J \in \D{N} \}$.

In order to translate Theorem \ref{thm:PrincipalClassification} to
Theorem \ref{thm:MeasurableClassification} we remark that
$J \subseteq K$ is equivalent to $(K) \subseteq (J)$, 
so that the inclusions in condition (1) have to be reversed.
The hypothesis $\Delta{N} \cap M = \emptyset$ can be reformulated 
as $\Delta{N} \subseteq I \minus M$, or equivalently 
$(I \minus M) \subseteq (\Delta{N})$, and replacing the principal filter 
$(\Delta{N})$ by the filter $\D{N}$ leads to the formulation 
of condition (2) in the general setting.

\begin{theorem} \label{thm:MeasurableClassification}
  Assume that $A$ contains three or more alternatives,
  and let $(I,\Sigma)$ be a measurable space as above.
  For every arrovian voting system $C \colon \PPO_A^{\smash{I,\Sigma}} \to \PPO_A$
  there exists a map $\d$ that associates to each set $N \in \Sigma$
  the filter $\d{N} \subseteq \Sigma^N$ of decisive sets relative to $N$.
  Moreover, $\D \colon \Sigma \to \FF{\Sigma}$, $\D{N} = \d{N} + N$, 
  is the unique map enjoying the following properties:
  \begin{enumerate}
  \item
    $\D{N} \subseteq (N)$ and $\D{N} \subseteq \D{M}$ 
    for all $M \subseteq N$ with $M,N \in \Sigma$.
    \hfill (Monotonicity)
  \item
    If $M \in \Sigma$ and $(I \minus M) \subseteq \D{N}$ 
    then $\D( N \sqcup M ) = \D{N} + M$.
    \hfill (Minimality)
  \item
    We have $a \wpref b$ if and only if $a \ndiff_N b$ 
    and $a \spref_{\d{N}} b$ for some $N \in \Sigma$.
    \hfill (Decisiveness)
  \end{enumerate}

  Conversely, given an arbitrary map $\D \colon \Sigma \to \FF{\Sigma}$ 
  satisfying (1), we set $\d{N} = \D{N} - N$ and define a map 
  $C_\D \colon \PPO_A^{\smash{I,\Sigma}} \to \PPO_A$ by the rule (3).
  The resulting voting system $C_\D$ satisfies 
  the arrovian axioms of unanimity, neutrality, monotonicity, 
  and independence of irrelevant alternatives.  Condition (2)
  ensures that $\d{N}$ is the filter of decisive sets relative to $N$.

  Finally, the voting system $C$ maps $\LPO_A^{\smash{I,\Sigma}}$ to $\LPO_A$ 
  if and only if for every $N \in \Sigma$ the filter $\d{N}$ 
  on $\Sigma^N$ is either trivial or an ultrafilter.
  \qed
\end{theorem}

Analogously to the finite case, Conditions (2) and (3) can be reformulated
as stated in Lemmas \ref{lem:WeakMinimality}--\ref{lem:WeakDecisiveness}.
The remaining verifications are a lengthy but routine transcription 
of the proof of Theorem \ref{thm:PrincipalClassification}, and will be omitted.

\begin{corollary}
  We have $C \subseteq C'$ if and only if $\D \subseteq \D'$.
  More explicitely, $C(P) \subseteq C'(P)$ for every profile 
  $P \in \PPO_A^{\smash{I,\Sigma}}$ if and only if
  $\D{N} \subseteq \D'{N}$ for every measurable set $N \in \Sigma$.
  \qed
\end{corollary}

  


The classification also holds for arrovian voting systems
$C \colon \LPO_A^{\smash{I,\Sigma}} \to \PPO_A$ defined on 
linear orderings, cf.\ Corollary \ref{cor:LinearClassification}.  
This leads again to Corollary \ref{cor:UniqueArrovianExtension},
which can be formulated without explicit reference to 
the coalition structure $\D$ nor any other technical details.





\bibliographystyle{plain}
\bibliography{junta}

\begin{thebibliography}{10}

\bibitem{Armstrong:1980}
T.~E. Armstrong.
\newblock Arrow's theorem with restricted coalition algebras.
\newblock {\em J. Math. Econom.}, 7(1):55--75, 1980.

\bibitem{Armstrong:1985}
T.~E. Armstrong.
\newblock Precisely dictatorial social welfare functions. {E}rratum and
  addendum to: ``{A}rrow's theorem with restricted coalition algebras''.
\newblock {\em J. Math. Econom.}, 14(1):57--59, 1985.

\bibitem{Arrow:1951}
K.~J. Arrow.
\newblock {\em Social {C}hoice and {I}ndividual {V}alues}.
\newblock Cowles Commission Monograph No. 12. John Wiley \& Sons Inc., New
  York, N. Y., 1951.
\newblock 2nd edition: John Wiley \& Sons Inc., 1963.

\bibitem{Bourbaki:1971}
N.~Bourbaki.
\newblock {\em \'{E}l\'ements de math\'ematique. {T}opologie g\'en\'erale.
  {C}hapitres 1 \`a 4}.
\newblock Hermann, Paris, 1971.

\bibitem{Fishburn:1970}
P.~C. Fishburn.
\newblock Arrow's impossibility theorem: concise proof and infinite voters.
\newblock {\em J. Econom. Theory}, 2:103--106, 1970.

\bibitem{Gibbard:1969}
A.~F. Gibbard.
\newblock Intransitive social indifference and the {A}rrow dilemma.
\newblock University of Chicago, 1969, unpublished manuscript.

\bibitem{Guilbaud:1952}
G.~Th. Guilbaud.
\newblock Les th\'eories de l'int\'er\^et g\'en\'eral et le probl\`eme logique
  de l'agr\'egation.
\newblock {\em \'Economie Appliqu\'ee}, 5:501--551, 1952.

\bibitem{KirmanSondermann:1972}
A.~P. Kirman and D.~Sondermann.
\newblock Arrow's theorem, many agents, and invisible dictators.
\newblock {\em J. Econom. Theory}, 5(2):267--277, 1972.

\bibitem{May:1952}
K.~O. May.
\newblock A set of independent necessary and sufficient conditions for simple
  majority decisions.
\newblock {\em Econometrika}, 20:680--684, 1952.

\bibitem{Sen:1986}
A.~Sen.
\newblock Social choice theory.
\newblock In {\em Handbook of mathematical economics, Vol.\ III}, volume~1 of
  {\em Handbooks in Econom.}, pages 1073--1181. North-Holland, Amsterdam, 1986.

\bibitem{Smith:1973}
J.~H. Smith.
\newblock Aggregation of preferences with variable electorate.
\newblock {\em Econometrica}, 41:1027--1041, 1973.

\bibitem{Weymark:1984}
J.~A. Weymark.
\newblock Arrow's theorem with social quasi-orderings.
\newblock {\em Public Choice}, 42:235--246, 1984.

\bibitem{Young:1974}
H.~P. Young.
\newblock An axiomatization of {B}orda's rule.
\newblock {\em J. Econom. Theory}, 9(1):43--52, 1974.

\end{thebibliography}

\end{document}